%
% Failure of flatness over finite-dimensional Hopf subalgebras
%	S. Skryabin
%
%	Plain TeX + AMS fonts ( amssym.def )
%
%	PARAMETERS
%
\hsize=5in
\baselineskip=12pt
\vsize=19.9cm
\parindent=10pt
\pretolerance=40
\predisplaypenalty=0
\displaywidowpenalty=0
\finalhyphendemerits=0
\hfuzz=2pt
\frenchspacing
\footline={\ifnum\pageno=1\else\hfil\tenrm\number\pageno\hfil\fi}
%
%	Additional Fonts
%
\input amssym.def
\font\titlefonts=cmbx12
\font\ninerm=cmr9
\font\ninebf=cmbx9
\font\ninei=cmmi9
\skewchar\ninei='177
\font\ninesy=cmsy9
\skewchar\ninesy='60
\font\nineit=cmti9
\def\reffonts{\baselineskip=0.9\baselineskip
	\textfont0=\ninerm
	\def\rm{\fam0\ninerm}%
	\textfont1=\ninei
	\textfont2=\ninesy
	\textfont\bffam=\ninebf
	\def\bf{\fam\bffam\ninebf}%
	\def\it{\nineit}%
        \rm}
%
%	Formatting
%
\def\frontmatter{\vbox{}\vskip1cm\bgroup
	\leftskip=0pt plus1fil\rightskip=0pt plus1fil
	\parindent=0pt
	\parfillskip=0pt
	\pretolerance=10000
	}
\def\endfrontmatter{\egroup\bigskip}
\def\title#1{{\titlefonts#1\par}}
\def\author#1{\bigskip#1\par}

\def\section#1\par{\ifdim\lastskip<\bigskipamount\removelastskip\fi
	\penalty-250\bigskip
	\vbox{\leftskip=0pt plus1fil\rightskip=0pt plus1fil
	\parindent=0pt
	\parfillskip=0pt
  \pretolerance=10000{\bf#1}}\nobreak\medskip
	}
\def\proclaim#1. {\medbreak\bgroup{\noindent\bf#1.}\ \it}
\def\endproclaim{\egroup
	\ifdim\lastskip<\medskipamount\removelastskip\medskip\fi}
\newdimen\itemsize
\def\setitemsize#1 {{\setbox0\hbox{#1\ }
	\global\itemsize=\wd0}}
\def\item#1 #2\par{\ifdim\lastskip<\smallskipamount\removelastskip\smallskip\fi
	{\leftskip=\itemsize
	\noindent\hskip-\leftskip
	\hbox to\leftskip{\hfil\rm#1\ }#2\par}\smallskip}
\def\Proof#1. {\ifdim\lastskip<\medskipamount\removelastskip\medskip\fi
	{\noindent\it Proof\if\space#1\space\else\ \fi#1.}\ }
\def\endproof{\hfill\hbox{}\quad\hbox{}\hfill\llap{$\square$}\medskip}
\def\Remark. {\ifdim\lastskip<\medskipamount\removelastskip\medskip\fi
        {\noindent\bf Remark. }}
\def\endremark{\medskip}
%
%       To LaTeX
%
\def\emph#1{{\it #1}\/}
\def\text#1{\hbox{#1}}
\def\mathrm#1{{\rm #1}}
\def\refeq#1{\/$(#1)$}
%
%		Citations
%
\newcount\citation
\newtoks\citetoks
\def\citedef#1\endcitedef{\citetoks={#1\endcitedef}}
\def\endcitedef#1\endcitedef{}
\def\citenum#1{\citation=0\def\curcite{#1}%
	\expandafter\checkendcite\the\citetoks}
\def\checkendcite#1{\ifx\endcitedef#1?\else
	\expandafter\lookcite\expandafter#1\fi}
\def\lookcite#1 {\advance\citation by1\def\auxcite{#1}%
	\ifx\auxcite\curcite\the\citation\expandafter\endcitedef\else
	\expandafter\checkendcite\fi}
\def\cite#1{\makecite#1,\cite}
\def\makecite#1,#2{[\citenum{#1}\ifx\cite#2]\else\expandafter\clearcite\expandafter#2\fi}
\def\clearcite#1,\cite{, #1]}
%
%		Reference section
%
\def\references{\section References\par
	\bgroup
	\parindent=0pt
	\reffonts
	\rm
	\frenchspacing
	\setbox0\hbox{99. }\leftskip=\wd0
	}
\def\endreferences{\egroup}
\newtoks\authtoks
\newif\iffirstauth
\def\checkendauth#1{\ifx\auth#1%
    \iffirstauth\the\authtoks
    \else{} and \the\authtoks\fi,%
  \else\iffirstauth\the\authtoks\firstauthfalse
    \else, \the\authtoks\fi
    \expandafter\nextauth\expandafter#1\fi
	}
\def\nextauth#1,#2;{\authtoks={#1 #2}\checkendauth}
\def\auth#1{\nextauth#1;\auth}
\newif\ifinbook
\newif\ifbookref
\def\nextref#1 {\par\hskip-\leftskip
	\hbox to\leftskip{\hfil\citenum{#1}.\ }%
	\initnextref}
\def\initnextref{\bookreffalse\inbookfalse\firstauthtrue\ignorespaces}
\def\paper#1{{\it#1},}
\def\InBook#1{\inbooktrue in ``#1",}
\def\book#1{\bookreftrue{\it#1},}
\def\journal#1{#1\ifinbook,\fi}
\def\BkSer#1{#1,}
\def\Vol#1{{\bf#1}}
\def\BkVol#1{Vol. #1,}
\def\publisher#1{#1,}
\def\Year#1{\ifbookref #1.\else\ifinbook #1,\else(#1)\fi\fi}
\def\Pages#1{\makepages#1.}
\long\def\makepages#1-#2.#3{\ifinbook pp. \fi#1--#2\ifx\par#3.\fi#3}
%
%	Math Stuff
%
\newsymbol\square 1003
\newsymbol\smallsetminus 2272
\let\setm\smallsetminus
\newsymbol\varnothing 203F
\newsymbol\vartriangleright 1342
\let\Rar\Rightarrow
\let\Lrar\Leftrightarrow
\let\ot\otimes
\let\sbs\subset
\let\sps\supset
\let\<\langle
\let\>\rangle
\def\co#1{^{\mkern1mu\mathrm{co}\mkern2mu#1}}
\def\coev{\mathrm{coev}}
\def\ev{\mathrm{ev}}
\def\id{\mathrm{id}}
\def\lco#1{\vphantom{H}\co{#1}\mkern-1mu}
\def\op{^{\mathrm{op}}}
\def\red{_{\mathrm{red}}}
\def\triv{_{\mathrm{triv}}}
\def\defop#1#2{\def#1{\mathop{\rm #2}\nolimits}}
\defop\chr{char}
\defop\Coend{Coend}
\defop\End{End}
\defop\Hom{Hom}
\defop\Img{Im}
\defop\Iso{Iso}
\defop\Ker{Ker}
\defop\pd{pd}
\def\limdir{\mathop{\vtop{\offinterlineskip\halign{##\hskip0pt\cr\rm lim\cr
	\noalign{\vskip1pt}
	$\scriptstyle\mathord-\mskip-10mu plus1fil
	\mathord-\mskip-10mu plus1fil
	\mathord\rightarrow$\cr}}}}
\def\liminv{\mathop{\vtop{\offinterlineskip\halign{##\hskip0pt\cr\rm lim\cr
	\noalign{\vskip1pt}
	$\scriptstyle\mathord\leftarrow
	\mskip-10mu plus1fil\mathord-\mskip-10mu plus1fil\mathord-$\cr}}}}
\def\mapr#1{{}\mathrel{\smash{\mathop{\longrightarrow}\limits^{#1}}}{}}
\def\lmapr#1#2{{}\mathrel{\smash{\mathop{\count0=#1
  \loop
    \ifnum\count0>0
    \advance\count0 by-1\smash{\mathord-}\mkern-4mu
  \repeat
  \mathord\rightarrow}\limits^{#2}}}{}}
\def\mapd#1#2{\llap{$\vcenter{\hbox{$\scriptstyle{#1}$}}$}\big\downarrow
  \rlap{$\vcenter{\hbox{$\scriptstyle{#2}$}}$}}
\def\diagram#1{\vbox{\halign{&\hfil$##$\hfil\cr #1}}}
\let\al\alpha
\let\be\beta
\let\ga\gamma
\let\de\delta
\let\ep\varepsilon
\let\ph\varphi
\let\si\sigma
\let\De\Delta
\def\frn{{\frak n}}
\def\frO{{\frak O}}
\def\calF{{\cal F}}
\def\calH{{\cal H}}
\def\calL{{\cal L}}
\def\calM{{\cal M}}
\def\bbA{{\Bbb A}}
\def\bbN{{\Bbb N}}
\def\bbZ{{\Bbb Z}}
\def\kov{{\overline k}}

\def\0{_{(0)}}
\def\1{_{(1)}}
\def\2{_{(2)}}
\def\3{_{(3)}}
\def\AM{\vphantom{\calM}_A\mskip1mu\calM}
\def\AMA{\AM^A}
\def\AMH{\AM^H}
\def\BMB{\vphantom{\calM}_B\mskip1mu\calM^B}
\def\MA{\calM_A}
\def\MAH{\calM_A^H}
\def\MB{\calM_B}
\def\MBH{\calM_B^H}

\citedef
Ag11
Ar-G03
Ber74
Br98
Br-G97
Br-Z10
Cae-G04
Dem-G
Doi83
Doi92
Doi-T86
Goo-W
Goo-Z10
Lar-Sw69
Ma94
Ma-W94
Mo
Nich-Z89a
Nich-Z89b
Ox61
Par72
Por11
Rot
Scha00
Schn92
Schn93
Sk06
Sk07
Sk08
Sk21
Sk-
Tak71
Tak72
Tak79
Wu-Z02
\endcitedef

\frontmatter
\title{Failure of flatness over finite-dimensional Hopf subalgebras}
\author{Serge Skryabin}
%\address{Institute of Mathematics and Mechanics,
%Kazan Federal University,\break
%Kremlevskaya St.~18, 420008 Kazan, Russia}
%\email{Serge.Skryabin@kpfu.ru}

\endfrontmatter

\bigskip
{\reffonts
\leftskip=20pt\rightskip=20pt
{\bf Abstract.}
It is proved in this paper that for any finite-dimensional nonsemisimple 
Hopf algebra $A$ there exists a Hopf algebra $H$ containing $A$ as a Hopf 
subalgebra such that $H$ is not flat over $A$. On the other hand, there is a 
class of infinite-dimensional Hopf algebras with the property that all Hopf 
algebras without exception are faithfully flat modules over Hopf subalgebras 
from this class.\par}

\section
Introduction

Let $H$ be a Hopf algebra over a field $k$. If $H$ is either commutative or 
cocommutative then $H$ is a faithfully flat module over any Hopf subalgebra. 
This is the earliest flatness result in the theory of Hopf algebras 
\cite{Tak72, Th. 3.1}. After its statement Takeuchi pointed out that it was an 
open problem whether the restriction on $H$ can be removed. The question was 
repeated two decades later in Montgomery's book \cite{Mo, Question 3.5.4}. 
However, Schauenburg showed that flatness does not hold in general 
\cite{Scha00}. His counterexamples use Hopf algebras whose antipodes are not 
bijective. Schauenburg proposed a modified question: {\it is every Hopf algebra 
with bijective antipode faithfully flat over all its Hopf subalgebras with 
bijective antipode}\/? As it turns out, this is also not true.

It is known that a Hopf algebra $H$ is a free module over its 
finite-dimensional Hopf subalgebras not only when $H$ is itself 
finite-dimensional \cite{Nich-Z89a}, but under the very mild restriction that 
$H$ is weakly finite \cite{Sk07}. However, in the present paper we will prove

\proclaim
Theorem 0.1.
For any nonsemisimple finite-dimensional Hopf algebra $A$ there exists a 
Hopf algebra $H$ containing $A$ as a Hopf subalgebra such that $H$ is neither 
left nor right flat over $A$.
\endproclaim

Here it is not essential whether the antipode of $H$ is bijective since we can 
always make a suitable replacement (see Lemma 1.2). Thus we do answer 
Schauenburg's question in the negative.

Whenever a Hopf algebra $H$ is a projective right module over its Hopf 
subalgebra $A$ it is known from a result of Doi \cite{Doi83} that all relative 
$(H,A)$-Hopf modules comprising the category $\MAH$ are projective $A$-modules. 
Thus, starting with a right $A$-module coalgebra $C$ which is not projective 
as an $A$-module, we can be sure that any Hopf algebra $H$ in which $A$ and 
$C$ are embedded is not right projective over $A$. If the antipode of $H$ is 
bijective, then such a Hopf algebra $H$ is even not flat over $A$ by the 
Masuoka-Wigner theorem \cite{Ma-W94}.

In section 2 we construct a Hopf algebra $H(C,A)$, universal with respect to 
the condition that there is a pair of compatible maps
$$
\imath:C\to H(C,A)\qquad\text{and}\qquad\jmath:A\to H(C,A)
$$
which are homomorphisms, respectively, of coalgebras and Hopf algebras. We call 
$H(C,A)$ the \emph{free Hopf algebra on a right $A$-module coalgebra $C$}. This 
generalizes Takeuchi's free Hopf algebra $H(C)$ on a coalgebra $C$ \cite{Tak71}. 

In section 3 it will be proved that the canonical maps $\imath$ and $\jmath$ 
are injective in the case when $\dim A<\infty$, but only for a special class 
of $A$-module coalgebras $C$. This is the most difficult part of the whole 
story. Injectivity can be recognized by means of algebra homomorphisms 
$H(C,A)\to B$. However, quite differently from the situation with the 
Hopf algebra $H(C)$, simple choices for the algebra $B$ now do not work. 
If $\dim B<\infty$, then any homomorphism $H(C,A)\to B$ factors through a 
residually finite-dimensional Hopf algebra, but such a Hopf algebra is 
faithfully flat over its finite-dimensional Hopf subalgebras by a result 
from \cite{Sk07}, and we do not get the desired conclusion. Therefore we have 
to work with the algebra $B=\End_kU$ where $U$ is an infinite-dimensional 
vector space. In other words, we verify existence of infinite-dimensional 
$H(C,A)$-modules on which $A$ acts faithfully.

A short section 1 collects all material which shows that the conclusions 
established in sections 2 and 3 lead to the proof of Theorem 0.1. An early 
indication that Theorem 0.1 can be true is traced to a work of Nichols and 
Zoeller \cite{Nich-Z89b}. There a Hopf algebra $H$ satisfying the universality 
property of $H(C,A)$ for some $A$ and $C$ of small dimension was defined by an 
explicit set of generators and relations. Since $A$ in their example was 
semisimple, the flatness over $A$ holds automatically. But even in this case 
one finds Hopf modules in the category $\MAH$ which are not free over $A$, and 
this means that something goes differently to what one would expect.

Theorem 0.1 may cast a pessimistic view on the flatness problem for an 
arbitrary Hopf algebra $H$. But things are really not bad when some 
infinite-dimensional Hopf subalgebras are considered. On the positive side we 
will prove

\proclaim
Theorem 0.2.
Let $A$ be a Hopf subalgebra of a Hopf algebra $H$. Suppose that $A$ is a 
semiprime Noetherian ring of finite global dimension. Then $H$ is right and 
left faithfully flat over $A$.
\endproclaim

Finiteness of global dimension in this theorem is a necessary assumption, as 
is seen from preceding Theorem 0.1. Necessity of the other assumptions about 
$A$ is less clear. In section 4 we will prove a more general theorem on 
projectivity of Hopf modules over $H$-comodule algebras. A longer list of 
ring-theoretic conditions used there weakens to some extent the conditions in 
Theorem 0.2. One may wonder whether the assumption that $A$ is semiprime 
Noetherian can be dropped altogether.

Brown and Goodearl asked whether arbitrary Noetherian Hopf algebras of finite 
global dimension are semiprime \cite{Br-G97, 1.9}. If $A$ has a quasi-Frobenius 
classical quotient ring $Q(A)$, then finiteness of its global dimension implies 
that $Q(A)$ is semisimple, and then $A$ is indeed semiprime. By \cite{Sk21, Th. 5.5} 
this argument applies to all residually finite-dimensional Noetherian Hopf 
algebras. For some time there stood the opposite question, asked by Wu and 
Zhang, whether every semiprime Noetherian Hopf algebra over a field of 
characteristic 0 has finite global dimension \cite{Wu-Z02, Question 0.4}. 
However, Goodearl and Zhang found counterexamples \cite{Goo-Z10, Remark 1.7} 
(see also \cite{Br-Z10, Example 7.3}).

A ring $R$ is said to be \emph{reduced} if $R$ has no nonzero nilpotent elements. 
A commutative algebra $A$ over a field $k$ is \emph{geometrically reduced} if 
the ring $A\ot_kK$ is reduced for every extension field $K$ of $k$.

\proclaim
Theorem 0.3.
Suppose that $A$ is a geometrically reduced commutative Hopf subalgebra of a 
Hopf algebra $H$. Then $H$ is right and left faithfully flat over $A$.
\endproclaim

Over a field $k$ of characteristic 0 every commutative Hopf algebra is known 
to be geometrically reduced. In this case any Hopf algebra $H$ is faithfully 
flat over all commutative Hopf subalgebras. This statement for $k$ the field 
of complex numbers appeared in a paper of Arkhipov and Gaitsgory \cite{Ar-G03, 
Prop. 3.12}. In our approach Theorem 0.3 is deduced almost immediately from 
Theorem 0.2.

If $\chr k>0$, then there exist nonsemisimple finite-dimensional commutative 
Hopf algebras. By Theorem 0.1 flatness over such Hopf subalgebras may fail. 
It will be nevertheless proved in section 5 that weakly finite Hopf algebras 
are faithfully flat over all commutative Hopf subalgebras. Recall that a ring 
$R$ is said to be \emph{weakly finite} or \emph{stably finite} if all 
surjective endomorphisms of any finitely generated free $R$-module are 
automorphisms. This property is satisfied under various finiteness conditions 
on a ring $R$. For example, residually finite-dimensional algebras over a 
field are weakly finite. Earlier it was proved that residually 
finite-dimensional Hopf algebras are faithfully flat over commutative Hopf 
subalgebras \cite{Sk08, Th. 0.3, Th. 7.13}. Thus we improve that result, 
applying different technique.

All work done in sections 4 and 5 is completely independent of the preceding 
part of the paper aimed at proving Theorem 0.1.

\section
1. The main result briefly explained

We assume that $k$ is the base field for algebras, coalgebras, and Hopf algebras. 
The subscript $k$ in the notation for the tensor product functor $\ot_k$ will 
be omitted.  Comultiplications will be denoted everywhere by the letter $\De$, 
counits by $\ep$, and antipodes by $S$. For a Hopf algebra $A$ its 
comultiplication $A\to A\ot A$ makes the categories $\AM$, $\MA$ of left and 
right $A$-modules into monoidal ones. A \emph{left $A$-module algebra} is an 
algebra in $\AM$. A \emph{right $A$-module coalgebra} is a coalgebra in $\MA$ 
(see \cite{Mo}).

If $\ph:A\to H$ is a homomorphism of Hopf algebras, then $H$ acquires a right 
$A$-module structure such that each element $a\in A$ acts on $H$ by means of 
the right multiplication by $\ph(a)$. With respect to this structure $H$ is a
right $A$-module coalgebra. In this situation the category of relative Hopf 
modules $\MAH$ makes sense. An object of this category is a right $A$-module 
$M$ equipped with a right $H$-comodule structure $M\to M\ot H$ given by an 
$A$-linear map with respect to the tensor product of right $A$-module 
structures on $M$ and $H$. An object $M\in\MAH$ is called \emph{$A$-finite} 
if $M$ is finitely generated as an $A$-module. As a special case of results 
due to Doi \cite{Doi83} we have

\setitemsize(a)
\proclaim
Lemma 1.1.
Given a homomorphism of Hopf algebras $\ph:A\to H,$ the following conditions 
are equivalent:

\item(a)
$H$ is a projective right $A$-module with respect to the action of $A$ arising 
from $\ph,$

\item(b)
all objects of $\MAH$ are projective in $\MA\,,$

\item(c)
all $A$-finite objects of $\MAH$ are projective in $\MA\,$.

\endproclaim

\Proof.
The implication (a)$\,\Rar\,$(b) follows immediately from \cite{Doi83, Th.~4, 
Cor.~1}. The converse (b)$\,\Rar\,$(a) is clear since $H\in\MAH$, and 
(b)$\,\Rar\,$(c) is trivial. Finally, the implication (c)$\,\Rar\,$(b) is 
deduced easily from the fact that all objects of $\MAH$ are locally 
$A$-finite (see \cite{Sk08, Lemmas 5.7, 7.1}).
\endproof

For any Hopf algebra $H$ there exist a Hopf algebra $\widehat H$ with bijective 
antipode and a Hopf algebra homomorphism $\ph:H\to\widehat H$ such that for every 
Hopf algebra $F$ with bijective antipode every homomorphism of Hopf algebras 
$H\to F$ factors through $\ph$ in a unique way. Schauenburg \cite{Scha00, 
Prop. 2.7} constructs $\widehat H$ as the colimit $\limdir H_i$ of the directed 
system of Hopf algebras
$$
H_0\lmapr2{f_0}H_1\lmapr2{f_1}H_2\lmapr2{}\cdots\eqno(1.1)
$$
where $H_i=H$ and $f_i=S^2$, the square of the antipode $S$ of $H$, for all 
$i\ge0$. The canonical homomorphism $\ph:H_0\to\limdir H_i$ satisfies the 
required universality property. It follows from this description that
$$
\Ker\ph=\bigcup_{n=1}^\infty\Ker S^n.
$$
If $A$ is a Hopf subalgebra of $H$, then its antipode is $S|_A$, and if 
$S|_A$ is injective, then so is $\ph|_A$. In this case $A$ is isomorphic to 
its image $\ph(A)$ in $\widehat H$, and so we may identify $A$ with this Hopf 
subalgebra of $\widehat H$.

\setitemsize(a)
\proclaim
Lemma 1.2.
Let $H$ be a Hopf algebra and $\widehat H$ the corresponding Hopf algebra with 
bijective antipode. Suppose that $A$ is a Hopf subalgebra of $H$ whose 
antipode is bijective. Then the following conditions are equivalent:

\item(a)
\ $H$ is right flat over $A,$

\item(b)
\ $H$ is right faithfully flat over $A,$

\item(c)
\ $H$ is a projective generator in $\MA,$

\item(d)
\ $\widehat H$ is right flat over $A$.

\noindent
Moreover, these conditions are equivalent to their left versions.
If $\,\dim A<\infty,$ then conditions {\rm(a), (b), (c), (d)} are also 
equivalent to the stronger condition

\item(e)
\ $H$ is a free right/left $A$-module.

\endproclaim

\Proof.
The implications (c)$\,\Rar\,$(b)$\,\Rar\,$(a) are obvious. If the antipode $S$ of $H$ 
is bijective, then (a)$\,\Lrar\,$(b)$\,\Lrar\,$(c) by the Masuoka-Wigner Theorem 
\cite{Ma-W94, Th. 2.1} and also (b)$\,\Lrar\,$(c) by \cite{Schn92, Cor. 1.8}.

In particular, (d) implies that $\widehat H$ is projective in $\MA$. In this case 
all objects of the category $\MA^{\widehat H}$ are projective in $\MA$ by Lemma 1.1. 
Since $H$ and $H/A$ may be viewed as objects of this category, they are then 
both projective in $\MA$, which implies also that $A$ is an $\MA$-direct 
summand of $H$. In other words, (d)$\,\Rar\,$(c).

Suppose that (a) holds. The composite homomorphism $f_{0i}:H_0\to H_i$ in the 
directed system (1.1) is $S^{2i}$. Let $f_{0i}^*\,H_i$ stand for the right 
$A$-module $H_i=H$ with respect to the pullback action of $A$ 
arising from $f_{0i}$. Then
$$
\widehat H\cong\limdir f_{0i}^*\,H_i\quad\text{in $\MA$}.
$$
Since $S|_A$ is bijective by the hypothesis, $f_{0i}|_A$ is a ring automorphism 
of $A$. It follows that $f_{0i}^*\,H_i$ is a flat $A$-module for each $i$. 
Since directed colimits of flat modules are flat (see, e.g., \cite{Rot, Prop. 
5.34}), so is $\widehat H$. Thus (a)$\,\Rar\,$(d).

Since the antipode of $\widehat H$ is a ring antiautomorphism which maps $A$ 
onto $A$, the right flatness of $\widehat H$ over $A$ is equivalent to the left 
flatness. Thus condition (d) is left-right symmetric. But the left versions of 
conditions (a), (b), (c), (d) are precisely these conditions for the Hopf 
algebras $A$, $H$, $\widehat H$ with their multiplications and 
comultiplications changed to the opposite ones. Hence the left versions of 
conditions (a), (b), (c), (d) are equivalent to each other.

If $\,\dim A<\infty$ then (c)$\,\Rar\,$(e) by a result of Schneider 
\cite{Schn93, Th. 2.4}.
\endproof

\proclaim
Corollary 1.3.
Let $A$ be a Hopf algebra and $C$ a right $A$-module coalgebra such that $C$ 
is not projective in $\MA$. Suppose that $H$ is a Hopf algebra, $f:C\to H$ 
a coalgebra homomorphism, and $\ph:A\to H$ a Hopf algebra homomorphism such that
$$
f(ca)=f(c)\,\ph(a)\qquad\text{for all $\,a\in A\,$ and $\,c\in C$.}\eqno(1.2)
$$
Then $H$ is not projective in $\MA$ with respect to the right action of $A$ on 
$H$ arising from $\ph$. If $\ph$ is injective and the antipode of $A$ is 
bijective, then $H$ is neither left nor right flat over its Hopf subalgebra 
$\ph(A)$.
\endproclaim

\Proof.
Let $\De:C\to C\ot C$ be the comultiplication. The composite map 
$$
(\id\ot f)\circ\De:C\to C\ot H
$$
gives a right $H$-comodule structure which makes $C$ an object of $\MAH$. 
Therefore the first conclusion follows immediately from Lemma 1.1. The second 
one is derived by applying Lemma 1.2.
\endproof

The functor of passage to dual vector spaces takes left $A$-modules to right 
ones and finite-dimensional left $A$-module algebras to right $A$-module 
coalgebras. This leads to the following construction.

Let $V$ be a finite-dimensional left $A$-module and $\,\tau:A\to\End_kV\,$ the 
corresponding representation of $A$. The algebra $\,\End_kV\,$ is a left 
$A$-module algebra with respect to the action of $A$ defined by the rule
$$
a\vartriangleright\be=\sum\,\tau(a\1)\circ\be\circ\tau(Sa\2),\qquad 
a\in A,\ \be\in\End_kV.\eqno(1.3)
$$
Thus $\,(a\vartriangleright\be)(v)=\sum\,a\1\be\bigl(S(a\2)v\bigr)\,$ for 
$v\in V$. Hence
$$
\Coend V=(\End_kV)^*\eqno(1.4)
$$
is a right $A$-module coalgebra.

The antipode of $A$ gives rise to a functor $\MA\to\AM$ by means of which the 
right $A$-module $V^*$ may be viewed as a left $A$-module. With respect to the 
monoidal structure on the category $\AM$ there is an isomorphism 
$\,\End_kV\cong V\ot V^*\,$ in $\AM$. Recall that projective left $A$-modules 
remain projective after tensoring with arbitrary left $A$-modules (see 
\cite{Br-G97, 1.2} or \cite{Par72, Lemma 5}), and note that $V$ is isomorphic 
to a direct summand of $\,V\ot V^*\ot V\,$ since the composite of 
two $A$-linear maps
$$
V\cong k\ot V\lmapr6{\coev\ot\id}V\ot V^*\ot V\lmapr6{\id\ot\ev}V\ot k\cong V
$$
is the identity map (here $\,\coev:k\to V\ot V^*\,$ and $\,\ev:V^*\ot V\to k\,$ 
are the coevaluation and the evaluation maps). It follows that $V$ is projective 
if and only if so is $V\ot V^*$.

Suppose that $\,\dim A<\infty$. Then $A$ is a Frobenius algebra by 
\cite{Lar-Sw69}, and therefore projectivity of left and right $A$-modules is 
preserved under duality. It follows that
$$
\eqalign{
\text{$\Coend V\,$ is projective in $\MA$}
&\quad\Lrar\quad\text{$\End_kV\,$ is projective in $\AM$}\cr
&\quad\Lrar\quad\text{$V\,$ is projective in $\AM$}.
}
$$

\medskip
{\bf Proof of Theorem 0.1.}
Take a nonprojective left $A$-module $V$ which contains a cyclic free direct 
summand. The right $A$-module coalgebra $C=\Coend V$ is then not projective in 
$\MA$. In section 2 we will describe the Hopf algebra $H(C,A)$, universal with 
respect to the property that there are homomorphisms of coalgebras $C\to H$ 
and Hopf algebras $A\to H$ satisfying (1.2). These two canonical maps into 
$H(C,A)$ are injective by Proposition 3.1. Hence $A$ is identified with a Hopf 
subalgebra of $H(C,A)$, and $H(C,A)$ is not flat over $A$ by Corollary 1.3.  
\endproof

The next corollary answers a question asked in \cite{Sk-}. As usual, for each 
subalgebra $A\sbs H$ we denote by $A^+$ its ideal of codimension 1 consisting 
of all elements $a\in A$ such that $\ep(a)=0$.

\proclaim
Corollary 1.4.
Let $A$ and $H$ be as in Theorem 0.1 with the antipode of $H$ being bijective.  
Then $H$ is left and right coflat over its factor coalgebra 
$\,\overline H=H/HA^+,$ but not faithfully coflat.
\endproclaim

\Proof.
Since $A$ is a Frobenius algebra, it is selfinjective, whence $A$ splits off 
in $H$ as an $\AM$-direct summand, and also as an $\MA$-direct summand. By 
\cite{Doi92, Prop. 3.2, 3.3} or \cite{Sk-, Prop. 2.1} $H$ is injective in the 
categories of left and right $\overline H$-comodules, which is equivalent 
to $H$ being left and right coflat over $\overline H$. Moreover, $A$ coincides 
with the subalgebra $\co{\overline H}H$ of $H$ consisting of all elements 
invariant with respect to the left coaction of $\overline H$ on $H$, e.g. by 
\cite{Sk-, Cor. 4.5}. If $H$ were left or right faithfully coflat over 
$\overline H$, then $H$ would have to be faithfully flat over $A$ by \cite{Tak79, Th.~2} 
and \cite{Ma94, Th. 1.11}.
\endproof

If $A$ is a semisimple Hopf algebra, then all $A$-modules are projective. Even 
in this case we obtain a Hopf algebra which exhibits deviation from the 
expected behaviour:

\proclaim
Corollary 1.5.
For any finite-dimensional Hopf algebra $A$ of dimension at least $2$ there 
exists a Hopf algebra $H$ containing a Hopf subalgebra isomorphic to $A$ and a 
finite-dimensional subcoalgebra $C$ such that $CA=C,$ but $C$ is not a free 
right $A$-module.
\endproclaim

\Proof.
Let $V=V_0\oplus V_1$ be the direct sum of a cyclic free left $A$-module $V_0$ 
and a nonzero trivial module $V_1$. Freeness is preserved under duality and 
tensoring with arbitrary $A$-modules. Hence in the direct sum decomposition
$$
V\ot V^*\cong(V_0\ot V_0^*)\oplus(V_0\ot V_1^*)\oplus(V_1\ot V_0^*)
\oplus(V_1\ot V_1^*)
$$
the first 3 summands are free $A$-modules, while the last one is nonzero 
with trivial action of $A$. This means that $\,\End_kV\,$ is not a free 
$A$-module, and the $A$-module coalgebra $\,C=\Coend V\,$ is consequently not 
a free $A$-module either. We may take $H=H(C,A)$ applying Proposition 3.1.
\endproof

The special case of Corollary 1.5 where $A$ is the 2-dimensional group algebra 
is known from a work of Nichols and Zoeller \cite{Nich-Z89b}. On the other 
hand, if $A$ is a finite-dimensional right coideal subalgebra of a weakly 
finite Hopf algebra $H$, then all objects of $\MAH$ are free $A$-modules by 
\cite{Sk07, Th. 6.1}.

\section
2. The free Hopf algebra on a Hopf module coalgebra

Let $A$ be a Hopf algebra and $C$ a right $A$-module coalgebra. Consider the 
category whose objects are the triples $(H,f,\ph)$ where $H$ is a Hopf algebra, 
$f:C\to H$ a coalgebra homomorphism, and $\ph:A\to H$ a Hopf algebra 
homomorphism satisfying (1.2). A morphism from $(H,f,\ph)$ to another triple 
$(H',f',\ph')$ in this category is a Hopf algebra homomorphism $\,t:H\to H'$ 
such that $f'=t\circ f$ and $\ph'=t\circ\ph$. 

This category of triples has the initial object, and we call the respective 
Hopf algebra $H(C,A)$ the \emph{free Hopf algebra on the $A$-module coalgebra 
$C$}. This Hopf algebra coincides with the free Hopf algebra $H(C)$ 
constructed by Takeuchi \cite{Tak71} in the case when $A=k$ is the trivial 
Hopf algebra. We will give an explicit description of $H(C,A)$ in terms of 
$H(C)$ and $A$.

Given an algebra $B$, the vector space $\Hom_k(C,B)$ of linear maps $C\to B$ 
has two algebra structures arising from the comultiplication $\De$ in $C$ and 
its opposite comultiplication $\De\op$. Thus the two multiplications of linear 
maps $f,g$ in $\Hom_k(C,B)$ are defined in Sweedler's notation by the rules
$$
(f*g)(c)=\sum\,f(c\1)\,g(c\2),\qquad
(f\times g)(c)=\sum\,f(c\2)\,g(c\1)
$$
where $c\in C$. They are called the \emph{convolution} and the \emph{twist 
convolution} \cite{Doi-T86} of $f$ and $g$. The map $c\mapsto\ep(c)1$ is the 
identity element with respect to both algebra structures. We may now speak 
about $*$-invertible and $\times$-invertible linear maps $C\to B$. This 
terminology will be equally used when $C$ is replaced by some other coalgebra 
or Hopf algebra.

\setitemsize(ii)
\proclaim
Lemma 2.1.
Let $f,g:C\to B$ be two linear maps, $\ph:A\to B$ an algebra homomorphism, $n$ 
an even integer, $m$ an odd integer. Consider the following conditions:
$$
\openup1\jot
\eqalignno{
f(ca)&=f(c)\,\ph(S^na)\qquad\text{for all $\,a\in A\,$ and $\,c\in C$}&(2.1)\cr
g(ca)&=\ph(S^ma)\,g(c)\qquad\text{for all $\,a\in A\,$ and $\,c\in C$}&(2.2)
}
$$

\item(i)
Suppose that $g$ is the $*$-inverse of $f$. If $n\ge0,$ then \refeq{2.1} is 
equivalent to \refeq{2.2} with $m=n+1$.

\item(ii)
Suppose that $g$ is the $\times$-inverse of $f$. If $n>0,$ then \refeq{2.1} is 
equivalent to \refeq{2.2} with $m=n-1$.

\noindent
If the antipode of $A$ is bijective then these conclusions are valid also for 
$n\le0$.
\endproclaim

\Proof.
Define linear maps $\,f_1,\,f_2,\,g_1,\,g_2:C\ot A\to B\,$ by the rules
$$
\openup1\jot
\vcenter{\halign{$#$&${}#$&\qquad\qquad$#$&${}#$\cr
f_1(c\ot a)&=f(ca),&f_2(c\ot a)&=f(c)\,\ph(S^na),\cr
g_1(c\ot a)&=g(ca),&g_2(c\ot a&)=\ph(S^ma)\,g(c)\cr
}}
$$
for $a\in A$ and $c\in C$. The $A$-module structure map $C\ot A\to C$ is a 
coalgebra homomorphism. The induced map $\,\Hom_k(C,B)\to\Hom_k(C\ot A,B)\,$ is 
therefore an algebra homomorphism with respect to either the convolution or 
the twist convolution multiplication. Since $f_1,\,g_1$ are the images of 
$f,\,g$ under this map, they are $*$-inverses of each other in (i) and 
$\times$-inverses of each other in (ii). On the other hand, $S^n$ and $S^m$ 
are $*$-inverses of each other when $m=n+1$ and $\times$-inverses of each 
other when $m=n-1$. Hence in the respective case $g_2$ is the $*$-inverse or 
the $\times$-inverse of $f_2$. It follows that $f_1=f_2$ if and only if 
$g_1=g_2$.
\endproof

Further on the notation $A_1\coprod A_2$ will be used for the coproduct in the 
category of (associative unital) algebras. This coproduct is also known as the 
free product of two rings amalgamating a common subring, the field $k$ in our 
case (see \cite{Ber74}). Takeuchi \cite{Tak71, Lemma 30} pointed out that the 
coproducts in the category of bialgebras and in the category of Hopf algebras 
are given by the coproducts in the category of algebras. Moreover, it was proved 
by Agore \cite{Ag11} and Porst \cite{Por11} that the categories of bialgebras 
and Hopf algebras are cocomplete, and the forgetful functors to the category 
of algebras preserve colimits. We will need only the simplest instance of 
those conclusions:

\setitemsize(iii)
\proclaim
Lemma 2.2.
Let $A_1$, $A_2$ be two bialgebras.

\item(i)
There is a unique comultiplication on the algebra $A_1\coprod A_2$ with 
respect to which the two canonical maps $A_i\to A_1\coprod A_2$ are bialgebra 
homomorphisms.

\item(ii)
$A_1\coprod A_2$ is the coproduct of $A_1$ and $A_2$ in the category of 
bialgebras.

\item(iii)
If $A_1$ and $A_2$ are Hopf algebras, then so is $A_1\coprod A_2,$ and if the 
antipodes of both $A_1$ and $A_2$ are bijective, then so is the antipode of 
$A_1\coprod A_2$.

\endproclaim

By \cite{Tak71, Cor. 9} the canonical map $C\to H(C)$ is injective. So we may 
identify $C$ with a subcoalgebra of $H(C)$. By Lemma 2.2 the coproduct 
$\,H(C)\coprod A\,$ in the category of algebras is a Hopf algebra. Let
$$
\textstyle
\imath_0:H(C)\to H(C)\coprod A\qquad\text{and}\qquad\jmath_0:A\to H(C)\coprod A
$$
be the canonical Hopf algebra homomorphisms. Put
$$
\textstyle H(C,A)=(H(C)\coprod A)/I
$$
where $I$ is the ideal of $\,H(C)\coprod A\,$ generated by
$$
\{\imath_0(ca)-\imath_0(c)\,\jmath_0(a)\mid\ c\in C,\ a\in A\}.
$$
Define $\imath:H(C)\to H(C,A)$ and $\jmath:A\to H(C,A)$ to be the composites of 
$\imath_0$ and $\jmath_0$ with the canonical map $\,H(C)\coprod A\to H(C,A)$.

\proclaim
Lemma 2.3.
The ideal $I$ is a Hopf ideal. Hence $H(C,A)$ is a Hopf algebra.
\endproclaim

\Proof.
The linear maps $\,\xi,\,\eta:C\ot A\to H(C)\coprod A\,$ defined by the rules
$$
\xi(c\ot a)=\imath_0(ca),\qquad\eta(c\ot a)=\imath_0(c)\,\jmath_0(a)
$$
for $c\in C$ and $a\in A$ are coalgebra homomorphisms. It follows that 
$\,\Img(\xi-\eta)\,$ is a coideal of the Hopf algebra $H(C)\coprod A$. Since 
$I$ is the ideal generated by this coideal, it is a biideal.

Now we have $\imath(ca)=\imath(c)\,\jmath(a)$ for all $c\in C$ and $a\in A$. Since 
$\imath$ is an algebra homomorphism, the linear map $\imath\circ S:H(C)\to H(C,A)$ 
is the $*$-inverse of $\imath$. Hence $(\imath\circ S)|_C$ is the $*$-inverse of 
$\imath|_C$. Applying Lemma 2.1 with $f=\imath|_C$, $\,g=(\imath\circ S)|_C$ 
and $\ph=\jmath$, we deduce that 
$\imath\bigl(S(ca)\bigr)=\jmath(Sa)\,\imath(Sc)$ for all $c\in C$ and $a\in A$. 
This means that
$$
S\bigl(\imath_0(ca)-\imath_0(c)\,\jmath_0(a)\bigr)
=\imath_0\bigl(S(ca)\bigr)-\jmath_0(Sa)\,\imath_0(Sc)\in I
$$
for all $c$ and $a$. Hence $\,S(I)\sbs I$.
\endproof

\proclaim
Proposition 2.4.
Given a Hopf algebra $H,$ a coalgebra homomorphism $f:C\to H,$ and a Hopf 
algebra homomorphism $\ph:A\to H$ satisfying \refeq{1.2}, there is a unique 
Hopf algebra homomorphism $\,t:H(C,A)\to H\,$ such that 
$\,f=t\circ\imath\,|_C\,$ and $\,\ph=t\circ\jmath$.
\endproclaim

\Proof.
Since $f$ extends to a Hopf algebra homomorphism $H(C)\to H$ by the universality 
property of $H(C)$, there is a unique Hopf algebra homomorphism 
$$
\textstyle t_0:H(C)\coprod A\to H
$$
such that $\,f=t_0\circ\imath_0\,|_C\,$ and $\,\ph=t_0\circ\jmath_0$. 
Condition (1.2) ensures that $t_0$ vanishes on the ideal $I$. So $t_0$ induces 
the required $t$.
\endproof

Let $B$ be an algebra. Following Takeuchi \cite {Tak71} denote by $L(C,B)$ the 
set of all infinite sequences $f_0,f_1,f_2,\ldots$ of linear maps $C\to B$ 
such that $f_i$ is the $*$-inverse of $f_{i-1}$ for odd integers $i>0$ and 
$f_i$ is the $\times$-inverse of $f_{i-1}$ for even $i>0$. Such a sequence is 
completely determined by its 0-term $f_0$.

A linear map $f:C\to B$ extends to an algebra homomorphism $H(C)\to B$ if and 
only if $f$ is the 0-term of some sequence in $L(C,B)$ \cite{Tak71, Prop. 4}. 
This entails a similar property of $H(C,A)$:

\proclaim
Proposition 2.5.
The algebra homomorphisms $H(C,A)\to B$ are in a bijective correspondence with 
the pairs $(f,\ph)$ where $f:C\to B$ is the $0$-term of some sequence in 
$L(C,B)$ and $\ph:A\to B$ is an algebra homomorphism satisfying \refeq{1.2}.
\endproclaim

\Proof.
The algebra homomorphisms $H(C,A)\to B$ are in a bijective correspondence with 
the pairs of algebra homomorphisms $H(C)\to B$ and $A\to B$ such that the 
induced homomorphisms $\,H(C)\coprod A\to B\,$ vanish on the ideal $I$.
\endproof

\proclaim
Corollary 2.6.
Let $k'$ be an extension field of $k$. Then the canonical $k'\!$-Hopf algebra 
homomorphism $\,H(C\ot k',\,A\ot k')\to H(C,A)\ot k'\,$ is an isomorphism.
\endproclaim

\Proof.
It follows from Proposition 2.5 that for each $k'$-algebra $B$ the induced map 
gives a bijection between the sets of $k'$-algebra homomorphisms 
$H(C,A)\ot k'\to B$ and $H(C\ot k',\,A\ot k')\to B$ 
(cf. \cite{Tak71, Cor. 8}).
\endproof

\proclaim
Lemma 2.7.
Let $\ph:A\to B$ be an algebra homomorphism, and let $(f_i)$ be a sequence in 
$L(C,B)$ whose $0$-term $f=f_0$ satisfies \refeq{1.2}. Then each $f_i$ satisfies 
either \refeq{2.1} with $n=i$ or \refeq{2.2} with $m=i,$ depending on whether 
$i$ is even or odd.
\endproclaim

\Proof.
This is deduced by induction on $i$ by applying Lemma 2.1.
\endproof

Now let $\,C=\Coend V$ and $B=\End_kU$ where $U$ and $V$ are two vector spaces 
such that $\dim V<\infty$. By (1.4) $\,C$ is a coalgebra with the 
comultiplication dual to the multiplication in the finite-dimensional algebra 
$\,\End_kV$. There is a canonical linear bijection
$$
\Hom_k(C,B)\cong B\ot C^*\cong\End_kU\ot\End_kV\cong\End_k(U\ot V).\eqno(2.3)
$$
It gives an isomorphism of the convolution algebra $\Hom_k(C,B)$ onto the 
endomorphism algebra $\End_k(U\ot V)$. However, the second multiplication 
on $\Hom_k(C,B)$ requires a different algebra on the right hand side.

There is an antiisomorphism of $\,\End_kV$ onto $\,\End_kV^*$ which assigns to 
a linear operator $\be$ on $V$ the dual operator $\be^*$ on $V^*$. It gives rise 
to an isomorphism of algebras
$$
\eqalign{
\Hom_k(C,B)\cong B\ot(C^*)\op&\cong\End_kU\ot(\End_kV)\op\cr
&\cong\End_kU\ot\End_kV^*\cong\End_k(U\ot V^*)
}\eqno(2.4)
$$
with respect to the twist convolution multiplication on $\Hom_k(C,B)$. 
Note that (2.4) is the composite of (2.3) and the linear bijection
$$
\End_k(U\ot V)\to\End_k(U\ot V^*),\qquad f\mapsto f^\flat,\eqno(2.5)
$$
such that $(\al\ot\be)^\flat=\al\ot\be^*$ for all $\al\in\End_kU$ and 
$\be\in\End_kV$. We call the linear operator $f^\flat\in\End_k(U\ot V^*)$ 
the \emph{halfdual} of the linear operator $f\in\End_k(U\ot V)$.

In a similar way we define $g^\flat\in\End_k(U\ot V)$ for each 
$g\in\End_k(U\ot V^*)$. Then clearly $f^{\flat\flat}=f$.

This definition of halfdual operators in our paper is made with respect to the 
symmetric monoidal category of vector spaces regardless of any $A$-module 
structures considered later.

\proclaim
Lemma 2.8.
Denote by $L'(V,U)$ the set of all infinite sequences $f_0,f_1,f_2,\ldots$ such 
that $f_i\in\End_k(U\ot V)$ when $i$ is even, $f_i\in\End_k(U\ot V^*)$ when 
$i$ is odd, and
$$
f_i^{\,\flat}=f_{i-1}^{\,-1}\quad\text{for all $\,i>0$}.
$$
There is a canonical bijection between $L(C,B)$ and $L'(V,U)$.
\endproclaim

\Proof.
This bijection is obtained by replacing each component $f_i$ of a sequence 
$(f_i)$ in $L(C,B)$ with its image, say $\tilde f_i$, under either (2.3) or 
(2.4), depending on whether $i$ is even or odd. By this assignment 
$\tilde f_i^{\,\flat}$ is the image of $f_i$ either under (2.3) when $i$ is 
odd or under (2.4) when $i$ is even. Since (2.3) and (2.4) are algebra 
isomorphisms with respect to the convolution and the twist convolution 
multiplications on $\Hom_k(C,B)$, the condition that $f_i$ is either 
$*$-inverse or $\times$-inverse of $f_{i-1}$ is equivalent to 
$\,\tilde f_i^{\,\flat}=\tilde f_{i-1}^{\,-1}$.
\endproof

Further on we assume that $U$ and $V$ are left $A$-modules. Let
$$
\ph:A\to\End_kU,\qquad\tau:A\to\End_kV
$$
be the algebra homomorphisms which give the action of $A$ on $U$ and $V$. The 
algebra $\End_kV$ is a left $A$-module algebra with respect to action (1.3). 
The dual spaces $V^*$ and $C=\Coend V$ are right $A$-modules in a natural way, 
and $C$ is moreover a right $A$-module coalgebra.

Denote by $U_{S^n}$ the $A$-module which has the same underlying vector space 
as $U$, but each element $a\in A$ acts on it by means of the operator 
$\ph(S^na)\in\End_kU$. It is either left or right $A$-module depending on 
whether $n$ is even or odd. Note that $U_{S^n}=U$ when $n=0$. Let $U\triv$ 
denote the vector space $U$ on which $A$ acts trivially, i.e., via the counit 
$\,\ep:A\to k$.

\proclaim
Lemma 2.9.
Under bijection \refeq{2.3} the linear maps $C\to B$ satisfying \refeq{2.1} 
correspond to the linear operators $f\in\End_k(U\ot V)$ such that
$$
\text{\rm$f\,$ is a morphism $\,U_{S^n}\ot V\to U\triv\ot V$ in $\AM$}.\eqno(2.6)
$$
Under bijection \refeq{2.4} the linear maps $C\to B$ satisfying \refeq{2.2} 
correspond to the linear operators $g\in\End_k(U\ot V^*)$ such that
$$
\text{\rm$g\,$ is a morphism $\,U\triv\ot V^*\to U_{S^m}\ot V^*$ in $\MA$}.\eqno(2.7)
$$
\endproclaim

\Proof.
Consider $\Hom_k(C,B)$ as an $A$-bimodule with respect to the left and right 
actions of $A$ defined by the rules
$$
(af)(c)=f(ca),\qquad(fa)(c)=f(c)\,\ph(a)
$$
where $a\in A$, $\,f\in\Hom_k(C,B)$, and $c\in C$. Identity (2.1) means 
precisely that
$$
af=fS^n(a)\quad\text{for all $a\in A$}.
$$
The corresponding two actions of $A$ on $\End_k(U\ot V)$ are as follows:
$$
\eqalign{
(af)(u\ot v)&=\sum\,\bigl(\id_U\ot\tau(a\1)\bigr)f\bigl(u\ot S(a\2)v\bigr),\qquad\cr
(fa)(u\ot v)&=f(au\ot v)
}
$$
where $a\in A$, $\,f\in\End_k(U\ot V)$, $\,u\in U$, and $v\in V$. 
The identity $af=fS^n(a)$ is expanded therefore as
$$
\sum\,\bigl(\id_U\ot\tau(a\1)\bigr)\,f\bigl(u\ot S(a\2)v\bigr)
=f\bigl(S^n(a)u\ot v\bigr).
$$
By means of standard transformations we obtain an equivalent identity
$$
\bigl(\id_U\ot\tau(a)\bigr)\,f(u\ot v)=\sum\,f\bigl(S^n(a\1)u\ot a\2v\bigr)
$$
which amounts to (2.6).

For an odd integer $m$ there is another bimodule structure on $\Hom_k(C,B)$ 
with the previously considered left action of $A$ but the new right action 
defined as
$$
(ga)(c)=\ph(S^ma)\,g(c),\qquad g\in\Hom_k(C,B).
$$
In terms of this structure identity (2.2) is rewritten as $ag=ga$.

Noting that for all $a\in A$, $\,\be\in\End_kV$, $\xi\in V^*$, and $v\in V$ we 
have
$$
\<(a\vartriangleright\be)^*(\xi),v\>
=\sum\,\<\xi,\,a\1\be\bigl(S(a\2)v\bigr)\>
=\sum\,\<\be^*(\xi a\1)S(a\2),\,v)\>,
$$
we deduce that $(a\vartriangleright\be)^*=a\vartriangleright\be^*$ where 
$a\vartriangleright\be^*\in\End_kV^*$ is defined by the rule
$$
(a\vartriangleright\be^*)(\xi)=\sum\,\be^*(\xi a\1)\,S(a\2).
$$
It follows that the two actions of $A$ on linear operators 
$g\in\End_k(U\ot V^*)$ obtained via (2.4) are given by the formulas
$$
\eqalign{
(ag)(u\ot\xi)&=\sum\,\bigl(\id_U\ot\tau'(Sa\2)\bigr)\,g(u\ot\xi a\1),\cr
(ga)(u\ot\xi)&=\bigl(\ph(S^ma)\ot\id_V\bigr)\,g(u\ot\xi)
}
$$
where $\tau':A\to\End_kV^*$ is the algebra antihomomorphism given by the 
right action of $A$ on $V^*$. The identity $ag=ga$ is equivalent to the 
identity
$$
g(u\ot\xi a)=\sum\,\bigl(\ph(S^ma\1)\ot\tau'(a\2)\bigr)\,g(u\ot\xi),
$$
i.e., to condition (2.7).
\endproof

We are ready to reformulate previously established facts concerning algebra 
homomorphisms $H(C,A)\to B$ in the presently considered case when $C=\Coend V$ 
and $B=\End_kU$. Recall that $\ph:A\to B$ stands for the representation of $A$ 
afforded by the $A$-module structure on $U$, and $\jmath:A\to H(C,A)$ is the 
canonical Hopf algebra map.

\proclaim
Lemma 2.10.
The algebra homomorphisms $\psi:H(C,A)\to B$ such that $\psi\circ\jmath=\ph$ are 
in a bijective correspondence with the linear operators $f\in\End_k(U\ot V)$ 
such that $f$ is the $0$-term of some sequence $(f_i)\in L'(V,U)$ and $f$ is a 
homomorphism of left $A$-modules $\,U\ot V\to U\triv\ot V$.

Moreover, each term $f_i$ of such a sequence satisfies either \refeq{2.6} with 
$n=i$ or \refeq{2.7} with $m=i,$ depending on whether $i$ is even or odd.  
\endproclaim

\Proof.
This lemma is a reformulation of Proposition 2.5 and Lemma 2.7 obtained by 
applying Lemmas 2.8 and 2.9. One can also prove directly, as in Lemma 3.4 later 
in the paper, that a linear operator $f\in\End_k(U\ot V)$ satisfies (2.6) if 
and only if $f^\flat$ is a morphism $U_{S^{n-1}}\ot V^*\to U\triv\ot V^*$ in 
$\MA$, while $g\in\End_k(U\ot V^*)$ satisfies (2.7) if and only if $g^\flat$ 
is a morphism $U\triv\ot V\to U_{S^{m-1}}\ot V$ in $\AM$.
\endproof

\proclaim
Lemma 2.11.
Let $\psi:H(C,A)\to B$ be the algebra homomorphism corresponding to some 
$f\in\End_k(U\ot V)$ as in Lemma 2.10. Suppose that the base field $k$ is 
algebraically closed. If $\,f(U\ot V')\not\sbs U\ot V'$ for each subspace $V'$ 
of $V$ other than the zero subspace and the whole $V,$ then the canonical map 
$\,\imath:C\to H(C,A)\,$ is injective.
\endproclaim

\Proof.
Since $\imath$ is a coalgebra homomorphism, its kernel $K=\Ker\imath$ is a 
coideal of $C$. The operator $f$ is the image of 
$\,\psi\circ\imath\in\Hom_k(C,B)\,$ under bijection (2.3). Since $K$ is 
contained in $\Ker(\psi\circ\imath)$, we get $f\in B\ot K^\perp$ where
$$
K^\perp=\{\be\in\End_kV\mid\<\be,K\>=0\}
$$
is a subalgebra of $\,\End_kV\cong C^*$. The condition on $f$ in the hypothesis 
means that $V$ is an irreducible $K^\perp$-module. This forces $K^\perp=\End_kV$ 
since $k$ is algebraically closed. Hence $K=0$.
\endproof

\Remark.
Suppose that $A$ is a Hopf algebra with bijective antipode. Let $\calH(C)$ and 
$\calH(C,A)$ be the Hopf algebras with bijective antipode obtained from $H(C)$ 
and $H(C,A)$ by applying Schauenburg's functor $H\mapsto\widehat H$. The Hopf 
algebra $\calH(C,A)$ satisfies the same universality property as $H(C,A)$, but 
with respect to pairs of homomorphisms from $C$ and $A$ to Hopf algebras with 
bijective antipode. It is the factor algebra of the coproduct $\calH(C)\coprod A$ 
by the ideal defined as in Lemma 2.3.

For each algebra $B$ denote by $\calL(C,B)$ the set of all sequences 
$(f_i)_{i\in\bbZ}$, infinite in both directions, of linear maps $C\to B$ such 
that $f_i$ is the $*$-inverse of $f_{i-1}$ for odd $i$ and the $\times$-inverse 
of $f_{i-1}$ for even $i$. Then the algebra homomorphisms $\calH(C,A)\to B$ 
are in a bijective correspondence with the pairs $(f,\ph)$ where $f:C\to B$ is 
the $0$-term of some sequence in $\calL(C,B)$ and $\ph:A\to B$ is an algebra 
homomorphism satisfying \refeq{1.2}. Under the same assumptions about $A$ and 
$C$ the proof of Proposition 3.1 can be extended to show that the canonical 
maps $C\to\calH(C,A)$ and $A\to\calH(C,A)$ are injective.
\endremark

\section
3. Injectivity of the canonical maps

The whole section is devoted to the proof of Proposition 3.1. With it we 
accomplish the crucial step in the proof of our main result stated in Theorem 0.1.

\proclaim
Proposition 3.1.
Let $A$ be a finite-dimensional Hopf algebra, and let $C=\Coend V$ where $V$ 
is a finite-dimensional left $A$-module containing an element with zero 
annihilator in $A$. Then the canonical maps $\imath:C\to H(C,A)$ and 
$\jmath:A\to H(C,A)$ are injective.
\endproclaim

By Corollary 2.6 the construction of the free Hopf algebra $H(C,A)$ commutes 
with base change. Therefore in the proof we may assume the base field $k$ to 
be uncountable and algebraically closed.

By the hypothesis $V$ contains a cyclic free submodule, say $V_0$. Since all 
finite-dimensional Hopf algebras are Frobenius, and therefore selfinjective, 
$V_0$ is a direct summand of $V$. Hence $V=V_0\oplus V_1$ where $V_1\sbs V$ is 
another submodule.

To show that $\jmath$ is injective we have to find an algebra homomorphism 
$$
\psi:H(C,A)\to B
$$
such that $\psi\circ\jmath$ is injective. We take $B=\End_kU$ where 
$U=U_0\oplus U_1$ is a left $A$-module written as the direct sum of a cyclic 
free submodule $U_0$ and a submodule $U_1$ on which $A$ acts trivially, i.e., 
via the counit $\ep:A\to k$.

By Lemma 2.10 $\psi$ is determined by a linear operator $f\in\End_k(U\ot V)$ 
such that $f=f_0$ is the $0$-term of some sequence $(f_n)\in L'(V,U)$ and $f$ 
intertwines the tensor product $A$-module structure on $U\ot V$ and another 
$A$-module structure given by the tensor product $U\triv\ot V$ in $\AM$. By 
the definition of $L'(V,U)$ the linear operator $f$ is invertible, whence $f$ 
gives an isomorphism $\,U\ot V\to U\triv\ot V$ in $\AM$.

If $V$ is not projective, then such an isomorphism cannot exist for a module 
$U$ with $U_1$ of finite dimension since in this case maximal projective direct 
summands of the $A$-modules $U\ot V$ and $U\triv\ot V$ have different dimension. 
If $U$ is infinite dimensional, then $\,U\ot V\cong U\triv\ot V$. However, the 
simplest choice of an isomorphism $f$ does not allow the construction of a 
sequence $(f_n)\in L'(V,U)$. Though it would suffice for our purposes to find 
only one sequence in $L'(V,U)$, we will not manage this as it appears to be 
difficult to control the explicit form of each term $f_n$. Instead, we will 
show that $\,L'(V,U)\ne\varnothing\,$ by applying topological arguments.

Let us take an $A$-module $U=U_0\oplus U_1$ as above with its trivial submodule 
$U_1$ of countably infinite dimension. The conclusion that $L'(V,U)$ is nonempty 
will be stated in Corollary 3.9 which sums up the work done in several preceding 
lemmas. Thus there is indeed at least one linear representation 
$\psi:H(C,A)\to\End_kU$ which extends the initially chosen representation of 
$A$ in $U$. Since $U$ is a faithful $A$-module, the map $\psi\circ\jmath\,$ is 
injective. Hence so is $\jmath$. This fact is already sufficient for the proof 
of Theorem 0.1.

The conclusion that $\imath$ is injective is based on the criterion stated in 
Lemma 2.11. Its verification will require extra work done at the end of this 
section.

\medskip
All topological considerations will refer to vector spaces of the form 
$T=\prod_{i\in I}T_i$ where $(T_i)_{i\in I}$ is an indexed collection of 
finite-dimensional vector spaces. We equip such a space $T$ with the 
\emph{pro-Zariski topology} which is the coarsest topology such that for each 
finite subset $F\sbs I$ the projection
$$
T\to T_F\quad\text{where $\,T_F=\prod_{i\in F}T_i$}
$$
is continuous when the finite-dimensional vector space $T_F$ is viewed with 
its Zariski topology (see Appendix at the end of the paper).

Here we need the assumption that the field $k$ is uncountable and algebraically 
closed. By Lemma A2 the space $T$ is then \emph{countably irreducible} with 
respect to the pro-Zariski topology, and so any countable intersection of 
nonempty open subsets is nonempty. Recall that \emph{$G_\de$-subsets} of a 
topological space are, by definition, countable intersections of open subsets. 
In a countably irreducible space any countable intersection of nonempty 
$G_\de$-subsets is itself a nonempty $G_\de$-subset.

Deriving existence properties from nonemptiness of $G_\de$-subsets is the 
essence of the Baire category method used many times in topology. If $k$ is 
the field of complex numbers then we could as well use the ordinary topology 
of finite-dimensional complex vector spaces and the product topology on the 
space $T$. This topology makes $T$ a Baire space since cartesian products of 
second countable Baire spaces are Baire by a theorem of Oxtoby 
\cite{Ox61, Th.~3}. Any countable intersection of dense $G_\de$-subsets in a 
Baire space is itself a dense $G_\de$-subset.

Suppose that $W=\bigoplus_{i\in I}W_i$ and $W'=\bigoplus_{i\in I}W'_i$ are two 
$A$-modules, either left modules or right ones, written as direct sums of 
collections of finite-dimensional submodules indexed by the same set $I$. The 
whole space $\Hom_A(W,W')$ does not have the topology we want. Its subspace
$$
\Hom_A^\oplus(W,W')
=\{f\in\Hom_A(W,W')\mid f(W_i)\sbs W'_i\text{ for all $i\in I$}\}
$$
is canonically isomorphic to the cartesian product 
$\prod_{i\in I}\Hom_A(W_i,W'_i)$ in which each component is finite dimensional, 
and so there is the pro-Zariski topology on the vector space 
$\Hom_A^\oplus(W,W')$, as introduced above.

There is some ambiguity in this notation as direct sum decompositions of $W$ 
and $W'$ should be understood from the context. Note that $\Hom_A^\oplus(W,W')$ 
is a subspace of $\,\Hom_k^\oplus(W,W')\cong\prod_{i\in I}\Hom_k(W_i,W'_i)\,$ 
with the induced topology.

\proclaim
Lemma 3.2.
Let $L$ be a vector subspace of finite codimension in $\,\Hom_A^\oplus(W,W')$. 
Then $L$ is closed in the pro-Zariski topology if and only if
$$
\{f\in\Hom_A^\oplus(W,W')\mid f(S)=0\}\sbs L
$$
for some finite-dimensional vector subspace $S$ of $\,W$.
\endproclaim

\Proof.
By Lemma A4 $L$ is closed if and only if there exists a finite subset 
$F\sbs I$ such that $L$ contains the kernel of the canonical map 
$\Hom_A^\oplus(W,W')\to\Hom_A(W_F,W'_F)$ where $W_F=\bigoplus_{i\in F}W_i$ and 
$W'_F=\bigoplus_{i\in F}W'_i$, i.e., $L$ contains all homomorphisms 
$f\in\Hom_A^\oplus(W,W')$ which vanish on $W_F$. Since each finite-dimensional 
vector subspace $S\sbs W$ is contained in $W_F$ for some $F$, this is 
equivalent to the condition in the statement of Lemma 3.2.
\endproof

\proclaim
Lemma 3.3.
Suppose that the index set $I$ is countable. Then the set\/ 
$\Iso_A^\oplus(W,W')$ of all isomorphisms in\/ $\Hom_A^\oplus(W,W')$ is a 
$G_\de$-subset of\/ $\Hom_A^\oplus(W,W')$ which is nonempty if and only if\/ 
$W_i\cong W'_i$ for each $i$. The assignment $f\mapsto f^{-1}$ gives a 
homeomorphism of\/ $\Iso_A^\oplus(W,W')$ onto\/ $\Iso_A^\oplus(W',W)$.  
\endproclaim

\Proof.
For each $i\in I$ the set $\Iso_A(W_i,W'_i)$ of all $A$-module isomorphisms 
$W_i\to W'_i$ is a Zariski open subset of the vector space $\Hom_A(W_i,W'_i)$. 
It follows that
$$
O_i=\{f\in\Hom_A^\oplus(W,W')\mid\quad f|_{\,W_i}\in\Iso_A(W_i,W'_i)\}
$$
is an open subset of $\Hom_A^\oplus(W,W')$. Clearly,
$\Iso_A^\oplus(W,W')=\bigcap_{i\in I}O_i\,$, and this subset of 
$\Hom_A^\oplus(W,W')$ with its induced topology is homeomorphic to the 
cartesian product $\prod_{i\in I}\Iso_A(W_i,W'_i)$ with the pro-Zariski 
topology. This explains the first assertion of the lemma.

There is a similar description of the set $\Iso_A^\oplus(W',W)$. Taking 
inverses gives an isomorphism of algebraic varieties 
$\,\Iso_A(W_i,W'_i)\cong\Iso_A(W'_i,W_i)\,$ for each $i$. Those isomorphisms 
are patched together to give a bicontinuous bijection of $\Iso_A^\oplus(W,W')$ 
onto $\Iso_A^\oplus(W',W)$.
\endproof

In our applications $W$ and $W'$ will be $A$-modules with the same underlying 
vector space, either $U\ot V$ or $U\ot V^*$, but different module structures. 
Making use of a countable basis of $U_1$ express $U_1$ as a direct sum of 
finite-dimensional subspaces
$$
\textstyle
U_1=(\bigoplus\limits_{i=1}^\infty U^0_i)
\oplus(\bigoplus\limits_{i=0}^\infty U^1_i)
$$
so that $\dim U^0_i=\dim U_0=\dim A$ for all $i>0$ and 
$\dim U^1_i=\dim V_1$ for all $i\ge0$. Also put $U^0_0=U_0$. Then
$$
\textstyle
U=U_0\oplus U_1=(\bigoplus\limits_{i=0}^\infty U^0_i)
\oplus(\bigoplus\limits_{i=0}^\infty U^1_i),
$$
a direct sum of submodules, with $A$ acting nontrivially only on $U^0_0$. Now
$$
\eqalign{
U\ot V&=\bigoplus_{i=0}^\infty\,
\bigl((U^0_i\ot V_0)\oplus(U^0_i\ot V_1)\oplus(U^1_i\ot V_0)
\oplus(U^1_i\ot V_1)\bigr),\cr
U\ot V^*&=\bigoplus_{i=0}^\infty\,
\bigl((U^0_i\ot V_0^*)\oplus(U^0_i\ot V_1^*)\oplus(U^1_i\ot V_0^*)
\oplus(U^1_i\ot V_1^*)\bigr)
}\eqno(3.1)
$$
are written as direct sums of $A$-submodules. Each of these direct summands 
remains an $A$-submodule when the original action of $A$ in $U$ is twisted by 
some power $S^n$ of the antipode $S$ of $A$ or when the original action is 
replaced by the trivial action.

The $\Hom_A^\oplus$ spaces further on in the text are formed with respect to 
direct sum decompositions obtained from (3.1) by joining some pairs of direct 
summands into a single summand. The precise decompositions to use are 
determined each time either by conditions (3.2) or by (3.3).

For each even $n\ge 0$ denote by $X_n$ the subset of 
$\,\Hom_A(U_{S^n}\ot V,\,U\triv\ot V)\,$ consisting of all homomorphisms 
$\,f:U_{S^n}\ot V\to U\triv\ot V\,$ such that
$$
\openup1\jot
\eqalign{
f\bigl((U^0_i\ot V_0)\oplus(U^1_i\ot V_1)\bigr)
&\sbs(U^0_i\ot V_0)\oplus(U^1_i\ot V_1)
\quad\text{for each $i\ge0$},\cr
f\bigl(U^0_0\ot V_1\bigr)&\sbs U^1_0\ot V_0,\cr
f\bigl((U^1_{i-1}\ot V_0)\oplus(U^0_i\ot V_1)\bigr)&\sbs 
(U^1_i\ot V_0)\oplus(U^0_{i-1}\ot V_1)\quad\text{for each $i>0$}.
}\eqno(3.2)
$$
For each odd $n>0$ denote by $X_n$ the subset of 
$\,\Hom_A(U\triv\ot V^*,\,U_{S^n}\ot V^*)\,$ consisting of all homomorphisms 
$\,g:U\triv\ot V^*\to U_{S^n}\ot V^*\,$ such that
$$
\openup1\jot
\eqalign{
g\bigl((U^0_i\ot V_0^*)\oplus(U^1_i\ot V_1^*)\bigr)
&\sbs(U^0_i\ot V_0^*)\oplus(U^1_i\ot V_1^*)
\quad\text{for each $i\ge0$},\cr
g\bigl(U^1_0\ot V_0^*\bigr)&\sbs U^0_0\ot V_1^*,\cr
g\bigl((U^1_i\ot V_0^*)\oplus(U^0_{i-1}\ot V_1^*)\bigr)&\sbs 
(U^1_{i-1}\ot V_0^*)\oplus(U^0_i\ot V_1^*)\quad\text{for each $i>0$}.
}\eqno(3.3)
$$
Thus $X_n$ is either $\Hom_A^\oplus(U_{S^n}\ot V,\,U\triv\ot V)$ or 
$\Hom_A^\oplus(U\triv\ot V^*,\,U_{S^n}\ot V^*)$ with respect to certain direct 
sum decompositions with collections of submodules indexed by the disjoint 
union $I=\bbN\coprod\bbN$ of two copies of the set of nonnegative integers 
$\bbN$.

\smallskip
Denote by $Y_n$ either $\,\Hom_A^\oplus(U\triv\ot V,\,U_{S^n}\ot V)\,$ or 
$\,\Hom_A^\oplus(U_{S^n}\ot V^*,\,U\triv\ot V^*)$, depending on the parity of 
$n$, with respect to the same direct sum decompositions as those used in the 
definition of $X_n$. For each $n$ both $X_n$ and $Y_n$ are vector subspaces of 
either $\,\End_k(U\ot V)\,$ or $\,\End_k(U\ot V^*)$, depending on the parity 
of $n$.

\proclaim
Lemma 3.4.
For each $n\ge0$ the assignment $f\mapsto f^\flat$ defines a bicontinuous 
isomorphism of topological vector spaces $\,Y_n\to X_{n+1}$.
\endproclaim

\Proof.
Recall linear bijection (2.5). The direct sum decomposition $V=V_0\oplus V_1$ 
gives rise to the block decomposition
$$
\openup1\jot
\displaylines{
\End_kV=E_{00}\oplus E_{01}\oplus E_{10}\oplus E_{11},\cr
E_{ij}=\{\be\in\End_kV\mid\be(V_j)\sbs V_i\ \text{ and }\ \be(V_{j'})=0\}
}
$$
where we put $j'=1-j$ for $j=0$ and $j=1$. Passage to dual operators yields
$$
\openup1\jot
\displaylines{
\End_kV^*=E_{00}^*\oplus E_{01}^*\oplus E_{10}^*\oplus E_{11}^*,\cr
E_{ij}^*=\{\be^*\in\End_kV^*\mid\be^*(V_i^*)\sbs V_j^*\ \text{ and }\ 
\be^*(V_{i'}^*)=0\}.
}
$$

\smallskip
Suppose first that $n$ is even. Then 
$\,Y_n=\widetilde Y\cap\Hom_A(U\triv\ot V,\,U_{S^n}\ot V)\,$
where $\widetilde Y$ is the subspace of $\End_k(U\ot V)$ consisting of all 
linear operators $f$ such that
$$
\openup1\jot
\eqalign{
f\bigl((U^0_i\ot V_0)\oplus(U^1_i\ot V_1)\bigr)
&\sbs(U^0_i\ot V_0)\oplus(U^1_i\ot V_1)
\quad\text{for each $i\ge0$},\cr
f\bigl(U^1_0\ot V_0\bigr)&\sbs U^0_0\ot V_1,\cr
f\bigl((U^1_i\ot V_0)\oplus(U^0_{i-1}\ot V_1)\bigr)&\sbs 
(U^1_{i-1}\ot V_0)\oplus(U^0_i\ot V_1)
\quad\text{for each $i>0$}
}
$$
(these inclusions are reversed version of (3.2)). Note that
$$
\openup1\jot
\eqalign{
\quad&\widetilde Y=(F_{00}\ot E_{00})\oplus(F_{01}\ot E_{01})
\oplus(F_{10}\ot E_{10})\oplus(F_{11}\ot E_{11})\qquad\text{where}\cr
\noalign{\smallskip}
F_{00}&=\{\al\in\End_kU\mid\al(U^0_i)\sbs U^0_i\ \,\forall i\ge0,
\ \ \al(U^1_0)=0,
\ \ \al(U^1_i)\sbs U^1_{i-1}\ \forall i>0\},\cr
F_{01}&=F_{10}=\{\al\in\End_kU\mid\al(U^0_i)\sbs U^1_i\text{ and }
\al(U^1_i)\sbs U^0_i\ \ \forall i\ge0\},\cr
F_{11}&=\{\al\in\End_kU\mid\al(U^0_i)\sbs U^0_{i+1}\text{ and }
\al(U^1_i)\sbs U^1_i\ \ \forall i\ge0\}.
}
$$
Similarly, 
$X_{n+1}=\widetilde X\cap\Hom_A(U\triv\ot V^*,\,U_{S^{n+1}}\ot V^*)\,$
where
$$
\widetilde X=(F_{00}\ot E_{00}^*)\oplus(F_{01}\ot E_{01}^*)
\oplus(F_{10}\ot E_{10}^*)\oplus(F_{11}\ot E_{11}^*)
$$
is the subspace of $\End_k(U\ot V^*)$ consisting of all linear operators 
satisfying (3.3). It is clear from the definition of halfdual operators that 
$\widetilde X$ is the image of $\widetilde Y$ under bijection (2.5).

For each $f\in\End_k(U\ot V)$ and $u\in U$ define linear maps 
$f_u:V\to U\ot V$ and $f^\flat_u:V^*\to U\ot V^*$ setting
$$
f_u(v)=f(u\ot v)\quad\text{for $v\in V$},\qquad 
f^\flat_u(\xi)=f^\flat(u\ot\xi)\quad\text{for $\xi\in V^*$}.
$$
Identifying $U\ot V^*$ with $\Hom_k(V,U)$ and $U\ot V$ with $\Hom_k(V^*,U)$ by 
means of canonical linear bijections, we have
$$
f^\flat_u(\xi)(v)=f_u(v)(\xi).\eqno(3.4)
$$
Indeed, it suffices to check this formula assuming that $f=\al\ot\be$ for some 
linear operators $\al\in\End_kU$ and $\be\in\End_kV$, but then
$$
f^\flat_u(\xi)(v)=\<\be^*(\xi),v\>\,\al(u)=\<\xi,\,\be(v)\>\,\al(u)=f_u(v)(\xi),
$$
as claimed.

A linear operator $f\in\End_k(U\ot V)$ is a morphism $U\triv\ot V\to U\ot V$ 
in $\AM$ if and only if $f_u$ is a morphism $V\to U\ot V$ in $\AM$ for each 
$u\in U$, which means that
$$
f_u(av)(\xi)=\sum a\1\cdot\bigl(f_u(v)(\xi a_2)\bigr)\quad
\text{for all $a\in A$, $v\in V$, $\xi\in V^*$}.\eqno(3.5)
$$
Similarly, $f^\flat$ is a morphism $U\triv\ot V^*\to U_S\ot V^*$ in $\MA$ if 
and only if $f^\flat_u$ is a morphism $V^*\to U_S\ot V^*$ in $\MA$ for each 
$u\in U$, which means that
$$
f^\flat_u(\xi a)(v)=\sum S(a\1)\cdot\bigl(f^\flat_u(\xi)(a_2v)\bigr)\quad
\text{for all $a\in A$, $v\in V$, $\xi\in V^*$}.\eqno(3.6)
$$
Note that identity (3.5) is equivalent to (3.6). If (3.5) holds, then
$$
\eqalign{
f^\flat_u(\xi a)(v)=f_u(v)(\xi a)
&=\sum S(a\1)a\2\cdot\bigl(f_u(v)(\xi a\3)\bigr)\cr
&=\sum S(a\1)\cdot\bigl(f_u(a\2v)(\xi)\bigr)
=\sum S(a\1)\cdot\bigl(f^\flat_u(\xi)(a_2v)\bigr).
}
$$
By means of a similar transformation (3.5) is deduced from (3.6).

Replacing in the preceding paragraph the $A$-module $U$ with $U_{S^n}$, we 
deduce that $f$ is a homomorphism of left $A$-modules 
$U\triv\ot V\to U_{S^n}\ot V$ if and only if $f^\flat$ is a homomorphism of 
right $A$-modules $U\triv\ot V^*\to U_{S^{n+1}}\ot V^*$. It follows that map 
(2.5) restricts to an isomorphism of vector spaces $\,Y_n\to X_{n+1}$.

As is seen from (3.4), $f_u=0$ for some $u\in U$ if and only if $f^\flat_u=0$. 
Hence $f$ vanishes on $U'\ot V$ where $U'$ is a vector subspace of $U$ if and 
only if $f^\flat$ vanishes on $U'\ot V^*$. Since every finite-dimensional 
subspace of $U\ot V$ (respectively, of $U\ot V^*$) is contained in a subspace 
of the form $U'\ot V$ (respectively, $U'\ot V^*$) for some finite-dimensional 
subspace $U'\sbs U$, it follows from Lemma 3.2 that the map $f\mapsto f^\flat$ 
induces a bijection between the sets of closed vector subspaces of finite 
codimension in $Y_n$ and in $X_{n+1}$. This allows us to apply Corollary A6 
which shows that the map $\flat:Y_n\to X_{n+1}$ is bicontinuous.

\smallskip
In the case when $n$ is odd we proceed as in the case of an even $n$ with 
appropriate changes in the definitions. Here
$$
\eqalign{
Y_n&=\widetilde Y\cap\Hom_A(U_{S^n}\ot V^*,\,U\triv\ot V^*),\cr
X_{n+1}&=\widetilde X\cap\Hom_A(U_{S^{n+1}}\ot V,\,U\triv\ot V)
}
$$
where $\widetilde X$ and $\widetilde Y$ are the subspaces, respectively, of 
$\End_k(U\ot V)$ and $\End_k(U\ot\nobreak V^*)$ consisting of all linear operators 
satisfying, respectively, (3.2) and the reversed version of (3.3). We have
$$
\openup1\jot
\eqalign{
\quad&\widetilde Y=(F_{00}\ot E_{00}^*)\oplus(F_{01}\ot E_{01}^*)
\oplus(F_{10}\ot E_{10}^*)\oplus(F_{11}\ot E_{11}^*)\cr
\quad&\widetilde X=(F_{00}\ot E_{00})\oplus(F_{01}\ot E_{01})
\oplus(F_{10}\ot E_{10})\oplus(F_{11}\ot E_{11})\qquad\text{where}\cr
\noalign{\smallskip}
F_{00}&=\{\al\in\End_kU\mid\al(U^0_i)\sbs U^0_i\text{ and }
\al(U^1_i)\sbs U^1_{i+1}\ \ \forall i\ge0\},\cr
F_{01}&=F_{10}=\{\al\in\End_kU\mid\al(U^0_i)\sbs U^1_i\text{ and }
\al(U^1_i)\sbs U^0_i\ \ \forall i\ge0\},\cr
F_{11}&=\{\al\in\End_kU\mid\al(U^0_i)\sbs U^0_{i-1}\ \,\forall i>0,
\ \ \al(U^0_0)=0,\ \,\al(U^1_i)\sbs U^1_i\ \ \forall i\ge0\}.
}
$$
Hence $\widetilde X$ is the image of $\widetilde Y$ under map (2.5). 
Continuity is verified as in the case of even $n$. It remains to prove that
a linear operator $f\in\End_k(U\ot V)$ is a morphism 
$U_{S^{n+1}}\ot V\to U\triv\ot V$ in $\AM$ if and only if 
$f^\flat$ is a morphism $U_{S^n}\ot V^*\to U\triv\ot V^*$ in $\MA$. For this 
we can again use identity (3.4). It suffices to consider the case $n=1$ and 
then replace $U$ with $U_{S^{n-1}}$. So let $n=1$. The required conditions on 
$f$ and $f^\flat$ are now expressed by means of the identities
$$
\eqalign{
\sum\,f_{S^2(a\1)u}(a\2v)(\xi)&=f_u(v)(\xi a),\cr
\sum\,f^\flat_{S(a\1)u}(\xi a\2)(v)&=f^\flat_u(\xi)(av).
}
$$
We claim that they are equivalent to each other. If the first one is 
satisfied, then
$$
\eqalign{
f^\flat_u(\xi)(av)=f_u(av)(\xi)&=\sum\,f_{S^2(a\2)S(a\1)u}(a\3v)(\xi)\cr
&=\sum\,f_{S(a\1)u}(v)(\xi a\2)=\sum\,f^\flat_{S(a\1)u}(\xi a\2)(v).
}
$$
The opposite implication is similar.
\endproof

\proclaim
Lemma 3.5.
The subsets $X_n(1)\sbs X_n$ and $Y_n(1)\sbs Y_n$ consisting of all $A$-module 
isomorphisms, respectively, in $X_n$ and $Y_n$ are nonempty $G_\de$-subsets. 
Moreover, the assignment $f\mapsto f^{-1}$ gives a homeomorphism of $X_n(1)$ 
onto $Y_n(1)$.  
\endproclaim

\Proof.
This lemma is a special case of Lemma 3.3. The fact that $X_n(1)$ and $Y_n(1)$ 
are nonempty follows from the existence of isomorphisms of left $A$-modules
$$
\openup1\jot
\displaylines{
(U^0_i)_{S^n}\ot V_0\cong(U^0_i)\triv\ot V_0
\quad\text{for each $i\ge0$},\cr
(U^0_0)_{S^n}\ot V_1\cong U^1_0\ot V_0,\cr
U^1_{i-1}\ot V_0\cong U^1_i\ot V_0\quad\text{and}\quad
U^0_i\ot V_1\cong(U^0_{i-1})\triv\ot V_1\quad\text{for each $i>0$}
}
$$
when $n$ is even and isomorphisms of right $A$-modules
$$
\openup1\jot
\displaylines{
(U^0_i)\triv\ot V_0^*\cong(U^0_i)_{S^n}\ot V_0^*
\quad\text{for each $i\ge0$},\cr
U^1_0\ot V_0^*\cong(U^0_0)_{S^n}\ot V_1^*,\cr
U^1_i\ot V_0^*\cong U^1_{i-1}\ot V_0^*\quad\text{and}\quad
(U^0_{i-1})\triv\ot V_1^*\cong U^0_i\ot V_1^*\quad\text{for each $i>0$}.
}
$$
when $n$ is odd. Recall that all components $U_i^0$, $U_i^1$ with the 
exception of $U^0_0$ are trivial $A$-modules. They remain trivial when the 
original action of $A$ is composed with $S^n$. 
The $A$-module $(U^0_0)_{S^n}$ is free of rank 1 for any $n$, and 
so too are $V_0$ and $V_0^*$. Therefore the $A$-modules
$$
(U^0_i)_{S^n}\ot V_0\,,\quad(U^0_i)\triv\ot V_0\,,\quad
(U^0_i)\triv\ot V_0^*,\quad(U^0_i)_{S^n}\ot V_0^*
$$
are all free of the same rank equal to $\,\dim A$, and the $A$-modules
$$
(U^0_0)_{S^n}\ot V_1\,,\quad U^1_i\ot V_0\,,\quad
U^1_i\ot V_0^*,\quad(U^0_0)_{S^n}\ot V_1^*
$$
are all free of the same rank equal to $\,\dim V_1$.
\endproof

\proclaim
Lemma 3.6.
For each $n\ge0$ the assignment $f\mapsto(f^{-1})^\flat$ defines a 
homeomorphism of $X_n(1)$ onto a $G_\de$-subset of $X_{n+1}$.
\endproclaim

\Proof.
This follows from Lemmas 3.4 and 3.5.
\endproof

Now we define recursively
$$
X_n(i)=\{f\in X_n(1)\mid(f^{-1})^\flat\in X_{n+1}(i-1)\}\quad\text{for $i>1$}.
$$
In this way we obtain a decreasing sequence of subsets 
$\,X_n\sps X_n(1)\sps X_n(2)\sps...$.

\proclaim
Lemma 3.7.
For each $i>0$ the set $X_n(i)$ is a nonempty $G_\de$-subset of $X_n$.
\endproclaim

\Proof.
We proceed by induction on $i$. Suppose that $i>1$ and $X_{n+1}(i-1)$ is a 
nonempty $G_\de$-subset of $X_{n+1}$. By Lemma 3.6 the set 
$$
T=\{(f^{-1})^\flat\mid f\in X_n(1)\}
$$
is a nonempty $G_\de$-subset of $X_{n+1}$. Then $T\cap X_{n+1}(i-1)$ is a 
$G_\de$-subset of $X_{n+1}$ as well. This subset is nonempty since $X_{n+1}$ 
is countably irreducible. The set $X_n(i)$ is the inverse image of 
$T\cap X_{n+1}(i-1)$ with respect to the homeomorphism of $X_n(1)$ onto $T$ 
given by Lemma 3.6. Hence $X_n(i)$ is a nonempty $G_\de$-subset of $X_n(1)$. 
Any such a subset is the intersection with $X_n(1)$ of some $G_\de$-subset of 
$X_n$. So $X_n(i)$ is the intersection of two $G_\de$-subsets of $X_n$, whence 
it is itself a $G_\de$-subset of $X_n$.
\endproof

\proclaim
Lemma 3.8.
The set $X(\infty)=\bigcap_{i=1}^\infty X_0(i)$ is a nonempty $G_\de$-subset 
of $X_0$. Moreover, $X(\infty)$ consists precisely of those $f\in X_0$ that 
occur as the $0$-terms of sequences in $L'(V,U)$.
\endproclaim

\Proof.
By Lemma 3.7 each $X_0(i)$ is a nonempty $G_\de$-subset of $X_0$. Hence 
$X(\infty)$ is a $G_\de$-subset of $X_0$ as well. It is nonempty since $X_0$ 
is countably irreducible.

For each $f\in X_0(i)$ there is a sequence $f_0,\ldots,f_i$ of length $i+1$ 
such that $f_0=f$ and
$$
f_j\in X_j(i-j),\qquad f_{j+1}=(f_j^{-1})^\flat\eqno(3.7)
$$
for all $j$ such that $0\le j<i$. Furthermore, $f_i\in X_i$ by Lemmas 3.4 and 
3.5. In order for this sequence to admit an extension by one term longer it is 
necessary and sufficient that $f_i$ be invertible, i.e., $f_i\in X_i(1)$, and 
then $f_j\in X_0(i+1-j)$ for $0\le j<i$.

Thus for $f\in X(\infty)$ there is an 
infinite sequence $(f_j)\in L'(V,U)$ which starts with $f_0=f$. Conversely, 
given a sequence $(f_j)\in L'(V,U)$ with $f_0\in X_0$, we get $f_0\in X_0(1)$ 
since $f_0$ is an invertible linear operator, and proceeding by induction we 
deduce that (3.7) holds for each $i>0$. Hence $f_0\in X(\infty)$.
\endproof

\proclaim
Corollary 3.9.
The set $L'(V,U)$ is nonempty.
\endproclaim

Each $f\in X(\infty)$ gives rise to an algebra homomorphism 
$\psi:H(C,A)\to\End_kU$. In general $\psi\circ\imath$ needn't be injective. We 
will show in Lemma 3.13 that each linear operator $f$ taken in a suitable 
nonempty open subset of $X(\infty)$ satisfies the hypothesis of Lemma 2.11. 
Let us start with the following observation:

\proclaim
Lemma 3.10.
There exists a nonempty open subset $O\sbs X_0$ such that for any $f\in O$ and 
any vector subspace $V'\sbs V$ with the property that
$$
f(U\ot V')\sbs U\ot V'\eqno(3.8)
$$
one has $V'=V'_0\oplus V'_1$ where $V'_i=V'\cap V_i,$ and 
moreover $V'=0$ whenever $V'_0=0,$ while $V'=V$ whenever $V'_0=V_0$.  
\endproclaim

\Proof.
We make use of the canonical surjections by means of which the space
$$
X_0=\Hom_A^\oplus(U\ot V,\,U\triv\ot V)
$$
is identified with a cartesian product of finite-dimensional vector spaces. 
All these maps are continuous by the definition of pro-Zariski topology on 
$X_0$. One such a map is
$$
\eqalign{
X_0\to&\Hom_A\bigl((U_0\ot V_0)\oplus(U^1_0\ot V_1),\,
({U_0}\,\triv\ot V_0)\oplus(U^1_0\ot V_1)\bigr)\cr
&\cong\Hom_A(U_0\ot V_0,\,{U_0}\,\triv\ot V_0)\oplus
\Hom_A(U_0\ot V_0,\,U^1_0\ot V_1)\cr
&\qquad\oplus\Hom_A(U^1_0\ot V_1,\,{U_0}\,\triv\ot V_0)\oplus
\Hom_A(U^1_0\ot V_1,\,U^1_0\ot V_1).
}
$$
Composing it with the projection onto the first direct summand we get 
a continuous surjective linear map
$$
\pi_1:X_0\to\Hom_A(U_0\ot V_0,\,{U_0}\,\triv\ot V_0).
$$
Another canonical surjection is
$$
\pi_2:X_0\to\Hom_A(U_0\ot V_1,\,U^1_0\ot V_0).
$$
And we also need a third canonical surjection
$$
\eqalign{
X_0\to&\Hom_A\bigl((U^1_0\ot V_0)\oplus(U^0_1\ot V_1),\,
(U^1_1\ot V_0)\oplus({U_0}\,\triv\ot V_1)\bigr)\cr
&\cong\Hom_A(U^1_0\ot V_0,\,U^1_1\ot V_0)\oplus
\Hom_A(U^1_0\ot V_0,\,{U_0}\,\triv\ot V_1)\cr
&\qquad\oplus\Hom_A(U^0_1\ot V_1,\,U^1_1\ot V_0)
\oplus\Hom_A(U^0_1\ot V_1,\,{U_0}\,\triv\ot V_1).
}
$$
Taking first the projection onto $\Hom_A(U^1_0\ot V_0,\,{U_0}\,\triv\ot V_1)$, 
consider the further composite map
$$
\pi_3:X_0\to\Hom_A(U^1_0\ot V_0,\,{U_0}\,\triv\ot V_1)
\to\Hom_A(U^1_0\ot V_0,\,V_1)
$$
where the second map is induced by the split epimorphism of $A$-modules 
$$
{U_0}\,\triv\ot V_1\to V_1
$$
arising from any fixed nonzero linear function $U_0\to k$. Thus $\pi_3$ is 
also surjective and continuous. The target spaces in the definition of 
$\pi_1,\pi_2,\pi_3$ are considered with the Zariski topology.

The spaces $\Hom_A(U_0\ot V_0,\,{U_0}\,\triv\ot V_0)$ and 
$\Hom_A(U_0\ot V_1,\,U^1_0\ot V_0)$ contain open subsets consisting of all 
isomorphisms in the respective spaces of homomorphisms. These subsets are 
nonempty since
$$
U_0\ot V_0\cong{U_0}\,\triv\ot V_0\quad\text{and}\quad 
U_0\ot V_1\cong U^1_0\ot V_0.
$$
In fact the $A$-modules in these isomorphisms are free of rank equal to 
$\dim U_0$ for the first pair of modules and $\dim V_1$ for the second pair. 

The space $\Hom_A(U^1_0\ot V_0,\,V_1)$ contains a nonempty open subset 
consisting of all surjective homomorphisms $U^1_0\ot V_0\to V_1$. 
We thus obtain 3 nonempty open subsets of $X_0\,$:
$$
\eqalign{
O_1&=\{f\in X_0\mid\pi_1(f)\text{ is an isomorphism}\},\cr
O_2&=\{f\in X_0\mid\pi_2(f)\text{ is an isomorphism}\},\cr
O_3&=\{f\in X_0\mid\pi_3(f)\text{ is surjective}\}.
}
$$

Suppose that $V'$ is a vector subspace of $V$ satisfying (3.8). Let 
$v=v_0+v_1$ where $v_0\in V_0$ and $v_1\in V_1$. For $u\in U_0=U^0_0$ we have
$$
f(u\ot v_0)\in(U_0\ot V_0)\oplus(U^1_0\ot V_1),\qquad
f(u\ot v_1)\in U^1_0\ot V_0
$$
by (3.2). Suppose that $v\in V'$. Then
$f(u\ot v)=f(u\ot v_0)+f(u\ot v_1)\in U\ot V'$ by (3.8), whence
$$
\eqalign{
f(u\ot v_0)+f(u\ot v_1)&\in(U\ot V')\cap\bigl((U_0\ot V_0)\oplus(U^1_0\ot V_1)
\oplus(U^1_0\ot V_0)\bigr)\cr
&\qquad=(U_0\ot V'_0)\oplus(U^1_0\ot V').
}
$$
It follows that $f(u\ot v_0)\in(U_0\ot V'_0)\oplus(U^1_0\ot V')$, and 
therefore
$$
\pi_1(f)(u\ot v_0)\in U_0\ot V'_0.
$$
This shows that $\pi_1(f)\bigl(U_0\ot p(V')\bigr)\sbs U_0\ot V'_0$ where 
$p:V\to V_0$ is the projection onto $V_0$ with kernel $V_1$.

Suppose that $f\in O_1$. Then $\pi_1(f)$ is an invertible linear operator on 
the vector space $U_0\ot V_0$, and we deduce that $\,\dim p(V')\le\dim V'_0$. 
However, this inequality yields $p(V')=V'_0$ since $V'_0=p(V'_0)\sbs p(V')$. 
In other words, $v_0\in V'_0$ for any vector $v=v_0+v_1$ lying in $V'$, but 
then also $v_1=v-v_0\in V'$, i.e., $\,V'=V'_0\oplus V'_1$.

If further $V'_0=0$, then $V'=V'_1$, and we get
$$
f(u\ot v)\in(U\ot V')\cap(U^1_0\ot V_0)=U^1_0\ot V'_0=0
$$
for all $u\in U_0$ and $v\in V'$. It follows that $V'=0$ when 
$f\in O_1\cap O_2$ since all endomorphisms in $O_2$ are injective on 
$U_0\ot V_1$.

Finally, for $u\in U^1_0$ and $v\in V'_0$ we have
$$
f(u\ot v)\in(U\ot V')\cap\bigl((U^1_1\ot V_0)\oplus(U_0\ot V_1)\bigr)
=(U^1_1\ot V'_0)\oplus(U_0\ot V'_1)
$$
by (3.2) and (3.8). Hence $\pi_3(f)(U^1_0\ot V'_0)\sbs V'_1$. If $V'_0=V_0$ 
and $f\in O_3$, then the previous inclusion entails $V'_1=V_1$, and so $V'=V$. 

We may take $O=O_1\cap O_2\cap O_3$. This is an open subset of $X_0$ 
satisfying all the required conditions. It is nonempty since $X_0$ is an 
irreducible topological space.
\endproof

\proclaim
Lemma 3.11.
Suppose that $L$ and $N$ are two finite-dimensional vector spaces of equal 
dimension. Then the set $O$ consisting of all linear operators 
$f\in\End_k(L\ot N)$ such that $f(L\ot N')\not\sbs L\ot N'$ for each 
nonzero proper vector subspace $N'$ of $N$ is a nonempty Zariski open 
subset of the vector space $\,\End_k(L\ot N)$.
\endproclaim

\Proof.
Let $n=\dim L=\dim N$. For each integer $d$ such that $0<d<n$ denote by 
$G_d(N)$ the Grassmann variety of $d$-dimensional vector subspaces of $N$. 
The set $Z_d$ consisting of all pairs $(f,N')$ where $f\in\End_k(L\ot N)$ 
and $N'\in G_d(N)$ satisfy
$$
f(L\ot N')\sbs L\ot N'
$$
is a closed subset of the product variety $\End_k(L\ot N)\times G_d(N)$. Since 
$G_d(N)$ is a projective variety, the projection 
$$
p:\End_k(L\ot N)\times G_d(N)\to\End_k(L\ot N)
$$
is a closed map. Hence $p(Z_d)$ is closed in $\End_k(L\ot N)$, and it follows 
that
$$
O=\End_k(L\ot N)\setm\bigcup_{0<d<n}p(Z_d)
$$
is indeed open in $\,\End_k(L\ot N)$. Taking the linear operator 
$f\in\End_k(L\ot N)$ defined by the formula
$$
f(u\ot v)=\si^{-1}(v)\ot\si(u),\qquad u\in L,\ v\in N,
$$
where $\si:L\to N$ is a linear bijection, we have 
$f(L\ot N')=\si^{-1}(N')\ot N$ for each vector subspace $N'$ of $N$, and it 
follows that $f\in O$. Thus $O$ is nonempty.

One can also note that
$$
\dim p(Z_d)\le\dim Z_d=n^4-n^2d(n-d)+d(n-d).
$$
In other words, $p(Z_d)$ is a closed subset of codimension at least 
$(n^2-1)d(n-d)$ in the whole space $\,\End_k(L\ot N)\,$ of dimension $n^4$.
\endproof

\proclaim
Lemma 3.12.
Suppose that $L$ and $N$ are two left $A$-modules, each being cyclic free. 
Then the set $O$ consisting of all $A$-linear maps $\,f:L\ot N\to L\triv\ot N\,$ 
such that
$$
f(L\ot N')\not\sbs L\ot N'
$$
for each nonzero proper vector subspace $N'$ of $N$ is a nonempty Zariski 
open subset of the vector space $\,\Hom_A(L\ot N,\,L\triv\ot N)$.
\endproclaim

\Proof.
Since $\,\Hom_A(L\ot N,\,L\triv\ot N)$ is a subspace of $\End_k(L\ot N)$, the 
openness of $O$ is immediate from Lemma 3.11. It remains to prove that $O$ is 
nonempty. We may assume that $L=N=A$ with the action of $A$ on itself by left 
multiplications. The $A$-module $A\ot A$ is generated freely by its subspace 
$1\ot A$. Hence there is an $A$-linear map $f:A\ot A\to A\triv\ot A$ such that 
$f(1\ot y)=y\ot1$ for all $y\in A$. Since
$$
x\ot y=\sum\,x\1\ot x\2S(x\3)y=\sum\,x\1\cdot\bigl(1\ot S(x\2)y\bigr)\qquad
\text{in $A\ot A$},
$$
we get
$$
f(x\ot y)=\sum\,S(x\2)y\ot x\1\qquad\text{for $x,y\in A$}.
$$
Suppose that $W$ is a vector subspace of $A$ such that $\,f(A\ot W)\sbs A\ot W$, 
i.e.,
$$
\sum\,S(x\2)y\ot x\1\in A\ot W\qquad\text{for all $x\in A$ and $y\in W$}.
$$
For $x=1$ this yields $W\ot1\sbs A\ot W$. Assuming that $W\ne0$, we get $1\in W$. 
But then we may take $y=1$ in the displayed containment above, and applying the 
map $\ep\ot\id$ where $\ep:A\to k$ is the counit of $A$, we deduce that 
$x=\sum\ep\bigl(S(x\2)\bigr)\,x\1$ lies in $W$ for each $x\in A$, i.e., $W=A$. 
Thus we have shown that $f\in O$. 
\endproof

\proclaim
Lemma 3.13.
There exists a nonempty open subset $\frO\sbs X_0$ such that
$$
f(U\ot V')\not\sbs U\ot V'
$$
for each $f\in\frO$ and each nonzero proper vector subspace $V'$ of $V$.
\endproclaim

\Proof.
We will use the continuous map $\pi_1$ defined in the proof of Lemma 3.10. By 
Lemma 3.12 applied with $L=U_0$ and $N=V_0$ there is a nonempty open subset 
$O_4\sbs X_0$ consisting of all $f\in X_0$ such that 
$$
\pi_1(f)(U_0\ot N')\not\sbs U_0\ot N'
$$
for each nonzero proper vector subspace $N'$ of $V_0$. Now take $\frO=O_4\cap O$ 
where $O$ is as in Lemma 3.10. The set $\frO$ is open in $X_0$ and is nonempty 
since $X_0$ is an irreducible topological space. 
Suppose that $V'$ is a vector subspace of $V$ satisfying (3.8) for some 
$f\in\frO$. Then
$$
\pi_1(f)(U_0\ot V'_0)\sbs U_0\ot V'_0,
$$
as we have seen in the proof of Lemma 3.10. But this inclusion implies that either 
$V'_0=0$ or $V'_0=V_0$. Hence either $V'=0$ or $V'=V$, again by Lemma 3.10.
\endproof

The proof of Proposition 3.1 is now completed as follows. Let 
$\imath:C\to H(C,A)$ be the canonical map. Since $X_0$ is countably irreducible, 
the set $X(\infty)\cap\frO$ is a nonempty $G_\de$-subset of $X_0$. By Lemma 3.8 
each $f\in X(\infty)\cap\frO$ is the 0-term of some sequence in $L'(V,U)$. 
Hence $f$ gives rise to an algebra homomorphism
$$
\psi:H(C,A)\to\End_kU,
$$
and by Lemma 2.11 $\imath$ is injective.

In retrospect our search for the terms of a sequence $(f_n)$ in the rather weird 
spaces $X_n$ does not seem to be unnatural. Starting from conditions (2.6) and 
(2.7) imposed by Lemma 2.10 we first made the choice that $f_n$ should map 
$U_0\ot V_1$ onto $U^1_0\ot V_0$ when $n$ is even and $U^1_0\ot V_0^*$ onto 
$U_0\ot V_1^*$ when $n$ is odd for some subspace $U^1_0\sbs U_1$. Then the 
other condition that $f_n=(f_{n+1}^{\,\flat})^{-1}$ forces $f_n$ to have many 
other nontrivial components.

Under the assumption that the base field $k$ is uncountable and algebraically 
closed we have shown that the Hopf algebra $H(C,A)$ has a left module of 
countably infinite dimension on which $A$ acts faithfully. In the case when $k$ 
is countable it is not clear whether every $A$-faithful $H(C,A)$-module must 
have uncountable dimension.

\section
4. Semiprime Noetherian comodule algebras and generalizations

Let $H$ be a Hopf algebra over the base field $k$. A \emph{right $H$-comodule 
algebra} is an algebra $A$ equipped with a right $H$-comodule structure 
$\rho_A:A\to A\ot H$ such that $\rho_A$ is an algebra homomorphism. For 
example, any right coideal subalgebra of $H$ is a right $H$-comodule algebra 
with respect to the comultiplication $\De$ in $H$. The main result of this 
section stated in Theorem 4.4 establishes projectivity of Hopf modules over a 
right $H$-comodule algebra $A$ under several assumptions.

For each right $H$-comodule algebra $A$ the categories of Hopf modules $\MAH$ 
and $\AMH$ are defined \cite{Doi92}. Both $A$ and $H$ are their objects with 
respect to the action of $A$ by either right or left multiplications. One of 
the assumptions required in the results of this section is that $A$ has no 
nonzero proper left ideals stable under the coaction of $H$, i.e., $A$ is a 
simple object of the category $\AMH$.

We will also need quotient rings. Recall that a ring $Q(A)$ containing $A$ as 
a subring is a \emph{classical left quotient ring} of $A$ if all nonzerodivisors 
of $A$ are invertible in $Q(A)$, and for each $x\in Q(A)$ there exists a 
nonzerodivisor $s\in A$ such that $sx\in A$. A left $A$-module $M$ is 
\emph{torsion} if
$$
Q(A)\ot_AM=0,
$$
which is equivalent to the condition that each element of $M$ is annihilated 
by a nonzerodivisor of $A$. The right version of $Q(A)$ is defined similarly.

These notions play a decisive role in several lemmas preceding Theorem 4.4.

\proclaim
Lemma 4.1.
Suppose that $A$ is a right $H$-comodule algebra such that $A$ is a simple 
object of $\AMH$. Then each nonzero $H$-subcomodule $U$ of a Hopf module 
$M\in\AMH$ has zero annihilator in $A$.

If $A$ has a classical left quotient ring, then $Q(A)\ot_AM\ne0$ unless $M=0$.
\endproclaim

\Proof.
The annihilator $I$ of $U$ in $A$ is stable under the coaction of $H$. Indeed,
$$
\sum a\0x\ot a\1=\sum a\0x\0\ot a\1x\1S(x\2)=\sum\rho_M(ax\0)\cdot S(x\1)=0
$$
for all $a\in I$ and $x\in U$. In this calculation $\rho_M:M\to M\ot H$ is the 
comodule structure map, and $M\ot H$ is regarded as a right $H$-module in a 
natural way. Thus $\rho_A(a)\in I\ot H$ for all $a\in I$, as claimed. Since 
$U\ne0$, we have $I\ne A$. In this case $I=0$ since $A$ is a simple 
object of $\AMH$.

Suppose that $M\ne0$. By the local finite-dimensionality of comodules there 
exists then a nonzero $H$-subcomodule $U\sbs M$ spanned linearly by finitely 
many elements. As we have proved, $U$ is not annihilated by any nonzero 
element of $A$. But each finite subset of a torsion left $A$-module is 
annihilated by some nonzerodivisor of $A$. Hence $M$ is not torsion.
\endproof

\proclaim
Lemma 4.2.
Let $A$ be a right $H$-comodule algebra. For each $A$-finite Hopf module 
$M\in\MAH$ the left $A$-module $\Hom_A(M,A)$ is an object of the category 
$\AMH$ with respect to a well-defined $H$-comodule structure, natural in $M$.
\endproclaim

\Proof.
Existence of comodule structures on $\Hom_A(M,A)$ was proved by Caenepeel and 
Gu\'ed\'enon \cite{Cae-G04, Prop. 4.2} for $A$-finite Hopf modules of the 
category $\AMH$. As pointed out in \cite{Sk-, Prop. 6.1}, the $\MAH$-version 
of that result does not require bijectivity of the antipode needed in 
\cite{Cae-G04}.

If $U$ is any finite-dimensional $H$-subcomodule of $M$ such that $M=UA$, then 
the right $A$-module $F=U\ot A$ is an object of $\MAH$ with respect to the 
tensor product of $H$-comodule structures on $U$ and $A$, and there is a 
canonical epimorphism $F\to M$ in $\MAH$. The dual space $U^*$ has an 
$H$-comodule structure which makes $U^*$ the left dual of $U$ in the monoidal 
category $\calM^H$ of right $H$-comodules. It follows that the left $A$-module 
$\Hom_A(F,A)\cong A\ot U^*$ is an object of $\AMH$ with respect to the 
tensor product of $H$-comodule structures on $A$ and $U^*$.

The left $A$-module $\Hom_A(M,A)$ embeds in $\Hom_A(F,A)$ as a submodule. By 
the proof of \cite{Sk-, Prop. 6.1} this submodule is stable also under the 
coaction of $H$. Hence $\Hom_A(M,A)$ is an $\AMH$-subobject of $A\ot U^*$.
\endproof

The classical left quotient ring $Q(A)$ is right flat over $A$. However, we 
need left flatness in the next lemma. For this reason we assume that $Q(A)$ 
satisfies all conditions required for the right classical quotient ring of $A$ 
as well as for the left one.

\setitemsize(2)
\proclaim
Lemma 4.3.
Let $A$ be a right $H$-comodule algebra such that $A$ is a simple object of 
$\AMH$ and $A$ has a von Neumann regular two-sided classical quotient ring $Q(A)$.

Let $N$ be an $\MAH$-subobject of a Hopf module $M\in\MAH$. If both $M$ and $N$ 
are finitely generated projective in $\MA,$ then $N$ is an $\MA$-direct summand 
of $M$.
\endproclaim

\Proof.
By Lemma 4.2 the canonical map
$$
\Hom_A(M,A)\to\Hom_A(N,A)\eqno(4.1)
$$
is a morphism of the category $\AMH$. Hence its cokernel $K$ is an 
object of $\AMH$ as well. We want to show that $K=0$.

Since $Q(A)$ is flat over $A$, there is an exact sequence of right $Q(A)$-modules
$$
0\to N\ot_AQ(A)\to M\ot_AQ(A)\to M/N\ot_AQ(A)\to0\eqno(4.2)
$$
in which the first two nonzero terms are finitely generated projective. Hence 
the third term is a finitely presented module. Since $Q(A)$ is von Neumann 
regular, all its modules are flat. Recall that all finitely presented flat 
modules are projective (see, e.g., \cite{Rot, Th. 3.56}). It follows that 
$M/N\ot_AQ(A)$ is projective, and therefore sequence (4.2) splits. Applying the 
functor $\,\Hom_{Q(A)}\bigl(?,\,Q(A)\bigr)$, we deduce that the canonical map
$$
\Hom_A\bigl(M,Q(A)\bigr)\to\Hom_A\bigl(N,Q(A)\bigr)
$$
is surjective. It follows that each right $A$-linear map $f:N\to A$ extends to 
a right $A$-linear map $g:M\to Q(A)$. The image $g(M)$ of the latter is a 
finitely generated right $A$-submodule of $Q(A)$. Since there exists then a 
nonzerodivisor $s\in A$ such that $s\,g(M)\sbs A$, this shows that the map 
$sf:N\to A$ extends to a right $A$-linear map $sg:M\to A$, i.e., $sf$ lies in 
the image of (4.1). In other words, $K$ is a torsion left $A$-module. By 
Lemma 4.1 $K=0$.

Thus (4.1) is a surjective map. It follows that the canonical map 
$$
\Hom_A(M,P)\to\Hom_A(N,P)
$$
is surjective for each finitely generated projective right $A$-module $P$. 
Taking $P=N$, we deduce that the inclusion $N\to M$ admits an $A$-linear 
splitting.
\endproof

A ring $R$ is called \emph{right coherent} if each finitely generated right 
ideal of $R$ is a finitely presented right $R$-module. For example, the 
polynomial algebra in infinitely many commuting indeterminates is coherent, but 
not Noetherian. In fact this algebra satisfies all the assumptions of the next 
theorem. Because of such examples we wish to relax Noetherian conditions on 
algebras. Further on we denote by $\,\pd_AV$ the \emph{projective dimension} 
of an $A$-module $V$. Basic properties of coherent rings and von Neumann 
regular rings, as well as the notions of global and weak dimensions are 
discussed, e.g., in \cite{Rot}.

\proclaim
Theorem 4.4.
Let $A$ be a right $H$-comodule algebra such that $A$ is a simple object of 
$\AMH$. Suppose that

\item(1)
$\,A$ has a von Neumann regular two-sided classical quotient ring $Q(A),$

\item(2)
$\,A$ is right coherent,

\item(3)
$\,\pd_AV<\infty$ for each finitely presented right $A$-module $V$.

\noindent
Then all objects of the category $\MAH$ are flat right $A$-modules. If there 
exists a ring homomorphism $A\to R$ for some nonzero right Artinian ring $R,$ 
then all objects of $\MAH$ are projective in $\MA$.
\endproclaim

\Proof.
We first prove that a Hopf module $M\in\MAH$ is projective in $\MA$ provided 
that $M$ is finitely presented in $\MA$. Put $d=\pd_AM$. Condition (3) ensures 
that $d<\infty$. We have to show that $d=0$.

If $U$ is any finite-dimensional $H$-subcomodule of $M$ such that $M=UA$, then 
the free $A$-module $F_0=U\ot A$ is an object of $\MAH$ such that the 
inclusion $U\to M$ extends to an epimorphism $\ph:F_0\to M$ in $\MAH$. Since 
$M$ is finitely presented, the kernel of $\ph$ is an $A$-finite subobject of 
$F_0$, and since $A$ is right coherent, $\Ker\ph$ is even finitely presented 
in $\MA$. By iterating we obtain an infinite exact sequence
$$
\cdots\to F_2\to F_1\to F_0\to M\to0
$$
in $\MAH$ such that each $F_i$ is a finitely generated free $A$-module. By the 
definition of projective dimension $d$ is the smallest integer such that the 
cokernel of the map $F_{d+1}\to F_d$ is projective in $\MA$. If $d>0$, then 
this cokernel is isomorphic to the image $N$ of $F_d$ in $F_{d-1}$. But then 
Lemma 4.3 shows that $N$ is an $\MA$-direct summand of $F_{d-1}$. Therefore the 
cokernel $F_{d-1}/N$ of the map $F_d\to F_{d-1}$ is projective in $\MA$ as well, 
yielding a contradiction. Thus our preliminary claim is proved.

Suppose that $M\in\MAH$ is finitely generated and projective in $\MA$. If $N$ is 
an $A$-finite $\MAH$-subobject of $M$, then the factor object $M/N$ is finitely 
presented in $\MA$, whence $M/N$ has to be projective in $\MA$, as we have just 
proved. This implies that $N$ is an $\MA$-direct summand of $M$, and therefore 
$N$ is finitely generated projective in $\MA$ too. For an arbitrary subobject 
$N\sbs M$ the set $\calF(N)$ of all $A$-finite subobjects of $N$ is directed 
by inclusion, and $N$ coincides with the union of those. It follows that 
$$
M/N\cong\limdir_{N'\in\calF(N)}M/N'
$$
is a directed colimit of projective $A$-modules, and therefore $M/N$ is flat 
in $\MA$. Thus all factor objects of $M$ are flat in $\MA$.

Each $A$-finite object of $\MAH$ is a factor object of some $F\in\MAH$ which 
is a finitely generated free $A$-module. Hence all $A$-finite objects of $\MAH$
are flat in $\MA$. Now an arbitrary object of $\MAH$, being the union of the 
directed set of its $A$-finite subobjects, is also flat in $\MA$.

One consequence of Lemma 4.2 is that the trace ideal
$$
\sum\nolimits_{f\in\Hom_A(M,A)}\,f(M)
$$
of an $A$-finite Hopf module $M\in\MAH$ is stable under the coaction of $H$ on 
$A$. If $M$ is nonzero and projective in $\MA$, then $\Hom_A(M,A)\ne0$, and 
the $H$-simplicity of $A$ implies that $M$ is a generator in $\MA$ 
\cite{Sk-, Cor. 6.4}. In this case $M\ot_AR\ne0$ for any ring homomorphism 
$A\to R$ provided that $R\ne0$.

Let $A\to R$ be a ring homomorphism where $R$ is a nonzero right Artinian 
ring. For each $A$-finite Hopf module $M\in\MAH$ the right $R$-module 
$M\ot_AR$ has finite length which we denote by $\ell(M)$. Clearly, $\ell(M)=0$ 
if and only if $M\ot_AR=0$. 

Suppose that $M$ is finitely presented in $\MA$. Since then $M$ is projective 
in $\MA$, we deduce that $\ell(M)>0$ whenever $M\ne0$. For any $A$-finite 
$\MAH$-subobject $N\sbs M$ the exact sequence $0\to N\to M\to M/N\to0$ splits 
in $\MA$. Hence exactness of the sequence is preserved after tensoring with 
$R$ over $A$, and therefore
$$
\ell(M)=\ell(N)+\ell(M/N).
$$
Since the Hopf factor module $M/N$ is finitely presented in $\MA$, it follows 
from this equality that $\ell(N)<\ell(M)$ whenever $N\ne M$. Given an 
ascending chain
$$
N_1\sbs N_2\sbs N_3\sbs\cdots\sbs M
$$
of $A$-finite subobjects of $M$, we have $\ell(N_i)\le\ell(N_{i+1})$ for all 
$i$, and $\ell(N_i)<\ell(N_{i+1})$ whenever $N_i\ne N_{i+1}$. Since 
$\ell(N_i)\le\ell(M)$ for all $i$, there can be only finitely many proper 
inclusions in this chain of subobjects. Hence $N_i=N_{i+1}$ for all 
sufficiently large $i$. Thus $M$ satisfies the ACC on $A$-finite subobjects. 
Since each object of $\MAH$ is a directed union of $A$-finite ones, it follows 
that all subobjects of $M$ are $A$-finite, and therefore all factor objects of 
$M$ are projective in $\MA$.

Each $A$-finite object of $\MAH$ is a factor object of some $\MA$-finitely 
presented one. Hence each $A$-finite object is projective in $\MA$, but 
then projectivity holds for all objects of the category $\MAH$ by 
\cite{Sk08, Lemmas 5.7, 7.1}.
\endproof

\Remark.
There is a version of Theorem 4.4 in which the condition that $A$ is a simple 
object of $\AMH$ is replaced by the weaker condition that $A$ is 
\emph{$H$-simple} in the sense that $A$ has no nonzero proper two-sided ideals 
stable under the coaction of $H$. However, we have to add another condition

\item(4)
$\,A$ has no faithful finitely generated torsion left $A$-modules

\noindent
to conditions (1), (2), (3) in the statement of Theorem 4.4. Then the 
conclusion of Lemma 4.3 will still be true. Indeed, the cokernel $K$ of map 
(4.1) is an $A$-finite object of $\AMH$ which is a torsion left $A$-module. 
By condition (4) $K$ is not faithful as an $A$-module. But each nonzero Hopf 
module of the category $\AMH$ is a faithful left $A$-module since its 
annihilator in $A$ is a two-sided ideal of $A$ stable under the coaction of 
$H$. It follows that $K=0$, and this allows us to repeat the proofs of Lemma 
4.3 and Theorem 4.4 without further changes.

Note that condition (4) is satisfied when $A$ is semiprime left bounded left 
Goldie (see \cite{Goo-W, Lemma 9.2}).
\endremark

\proclaim
Corollary 4.5.
Let $A$ be a right $H$-comodule algebra such that $A$ is a simple object of 
$\AMH$. All Hopf modules of the category $\MAH$ are projective in $\MA$ under 
any of the following additional assumptions:

\item(a)
$A$ is semiprime two-sided Noetherian of finite global dimension,

\item(b)
$A$ is a semiprime two-sided Goldie right coherent ring of finite weak dimension,

\item(c)
$A$ is von Neumann regular with at least one simple Artinian factor ring.

\endproclaim

\Proof.
All conditions of Theorem 4.4 are satisfied in each of cases (a), (b), (c). 
Case (a) is a subcase of (b). In case (b) $A$ has a semisimple Artinian 
classical quotient ring $Q(A)$ by the Goldie theorem \cite{Goo-W, Th. 6.15}. 
Thus condition (1) is satisfied. Since condition (2) is assumed, the projective 
dimension of any finitely presented right $A$-module $V$ coincides with its 
flat dimension, whence (3) holds by the finiteness of the weak dimension of $A$.

In case (c) $Q(A)\cong A$. Each finitely generated one-sided ideal of a von 
Neumann regular ring is generated by an idempotent, and each finitely 
presented module is projective. Hence conditions (2) and (3) also hold.  
\endproof

\proclaim
Corollary 4.6.
Let $A\sbs B\sbs H$ where $H$ is a Hopf algebra, $B$ a Hopf subalgebra with 
bijective antipode, and $A$ a right coideal subalgebra such that $A$ is an 
$\MA$-direct summand of $B$. Suppose that $A$ satisfies conditions 
$(1)$, $(2)$, $(3)$ in the statement of Theorem 4.4. Then all Hopf modules of 
the category $\MAH$ are projective in $\MA$. As a consequence, $H$ is right 
faithfully flat over $A$.
\endproclaim

\Proof.
The counit $\ep:H\to k$ gives by restriction a ring homomorphisms $A\to k$. Thus 
we have only to check that $A$ is a simple object of $\AMH$. Since $A\sbs B$, 
each $\AMH$-subobject $I$ of $A$ is a left ideal of $A$ such that 
$\De(I)\sbs I\ot B$. Hence the left ideal $BI$ of $B$ is a $\BMB$-subobject of 
$B$. Since the antipode of $B$ is bijective, the bialgebra $B\op$ obtained 
from $B$ by changing its multiplication to the opposite one is a Hopf algebra.
The fundamental theorem on Hopf modules \cite{Mo, Th. 1.9.4} applied to $B\op$ 
shows that $B$ is a simple object of $\BMB$. It follows that $BI$ is equal to 
either $B$ or 0. But $I=A\cap BI$ by the condition that $A$ is an $\MA$-direct 
summand of $B$. Hence either $I=A$ or $I=0$.

So Theorem 4.4 does apply. Since $H$ is an object of $\MAH$ and $A$ is its 
subobject, the factor object $H/A$ is projective in $\MA$. Hence $H$ is 
a projective generator in $\MA$, yielding faithful flatness.
\endproof

\proclaim
Corollary 4.7.
Let $A\sbs B\sbs H$ where $H$ is a Hopf algebra, $B$ a Hopf subalgebra with 
bijective antipode, $A$ a right coideal subalgebra such that $A$ is an 
$\AM$-direct summand of $B$. Suppose that

\item(1)
$A$ has a von Neumann regular two-sided classical quotient ring $Q(A),$

\item(2)
$A$ is left coherent,

\item(3)
$\pd_AV<\infty$ for each finitely presented left $A$-module $V,$

\noindent
Then all Hopf modules of the category $\AMH$ are projective in $\AM$. 
As a consequence, $H$ is left faithfully flat over $A$.
\endproclaim

\Proof.
With $\AMH$ and $\MAH$ interchanged bijectivity of the antipode of $H$ becomes 
an issue in the respective versions of Lemmas 4.1 and 4.2. We remove this 
issue by passing to the Hopf algebra with bijective antipode $\widehat H$ 
obtained by Schauenburg's construction. Recall that the kernel of the 
canonical Hopf algebra homomorphism $H\to\widehat H$ consists of all elements 
of $H$ annihilated by some power of the antipode $S$ of $H$. Since $S|_B$ is 
bijective, the Hopf subalgebra $B$ embeds in $\widehat H$. Hence so too does 
the right coideal subalgebra $A$. Each object of the category $\AMH$ may be 
regarded as an object of $\AM^{\widehat H}$. Replacing $H$ with $\widehat H$ 
we reduce the proof to the case when $S$ is bijective. Changing the 
multiplication to the opposite one we obtain then a Hopf algebra $H\op$ with 
its right coideal subalgebra $A\op$ satisfying all conditions of Theorem 4.4.  
The category $\AMH$ may be identified with $\calM_{A\op}^{H\op}$, and thus we 
are in the situation of Corollary 4.6.
\endproof

\medskip
{\bf Proof of Theorem 0.2.}
The fact that $A$ has an Artinian right or left classical quotient ring implies 
that the antipode of $A$ is bijective \cite{Sk06, Th.~A}. Therefore the 
conclusion of Theorem 0.2 is obtained as a special case of Corollaries 
4.6 and 4.7 in which we take $B=A$. All conditions of Theorem 4.4 are 
satisfied by case (a) of Corollary 4.5. Bijectivity of the antipode implies 
also that one-sided conditions on $A$ are equivalent to two-sided conditions.
\endproof

\section
5. Flatness over commutative Hopf subalgebras

This section provides extra arguments needed to deal with the question of 
flatness over arbitrary commutative Hopf subalgebras in positive characteristic. 
The final result will be stated in Theorem 5.4. However, Theorem 0.3 follows 
instantaneously from the already proved Theorem 0.2, as we can show 
straightaway.

\medskip
{\bf Proof of Theorem 0.3.}
Since each commutative Hopf algebra is a directed union of finitely generated 
Hopf subalgebras, the proof reduces to the case when $A$ is finitely generated. 
The assumption that $A$ is geometrically reduced implies that the affine group 
scheme $G$ represented by $A$ is smooth \cite{Dem-G, Ch. II, \S~5, Th. 2.1}. 
This means that all local rings of $A$ are regular. Hence $A$ has finite 
global dimension equal to the dimension of $G$, and it follows that $A$ 
satisfies all the assumptions of Theorem 0.2.
\endproof

\Remark.
Suppose that $A$ is a commutative Hopf subalgebra of a Hopf algebra $H$. 
If $\ga:A\to k$ is an algebra homomorphism and $\xi_\ga:H\to k$ is a linear 
map vanishing on the right ideal $(\Ker\ga)H$ of $H$, then one can consider 
the composite map
$$
H\lmapr2\De H\ot H\lmapr4{\xi_\ga\ot\id}k\ot H\cong H.
$$
In fact the action of $\xi_\ga\in H^*$ with respect to the standard right 
$H^*$-module structure on $H$ is defined precisely by this map. In the proof 
proposed by Arkhipov and Gaitsgory \cite{Ar-G03, Lemma 3.14} such a map was 
used to obtain an isomorphism between the original $A$-module structure on $H$ 
and the module structure twisted by $\ga$. However, bijectivity of the map 
$H\to H$ given by the action of $\xi_\ga$ is not clear unless $\xi_\ga$ is an 
invertible element of the dual algebra $H^*$.
\endremark

Now we are going to consider arbitrary commutative Hopf subalgebras. In doing 
this we will need a suitable extension of several arguments used in section 4. 
First we show that under the assumptions of the next lemma splitting of direct 
summands does not require the condition that $Q(A)$ is von Neumann regular in 
the hypothesis of Lemma 4.3.

\proclaim
Lemma 5.1.
Let $A$ be a finitely generated commutative Hopf subalgebra of a Hopf algebra 
$H,$ and let $N$ be an $\MAH$-subobject of an $A$-finite Hopf module 
$M\in\MAH$. If $N$ is finitely generated projective in $\MA,$ then $N$ is an 
$\MA$-direct summand of $M$.
\endproclaim

\Proof.
Finitely generated commutative Hopf algebras are Gorenstein \cite{Br98, Prop. 
2.3, Step 1}. Hence so is $A$, and therefore $A$ has a selfinjective Artinian 
classical quotient ring $Q(A)$. The projective $Q(A)$-module $N\ot_AQ(A)$ is 
then also injective, whence sequence (4.2) splits. Since the antipodes of 
commutative Hopf algebras are involutory, $A$ is a simple object of $\AMA$, 
and therefore also a simple object of $\AMH$. We can continue as in the proof 
of Lemma 4.3.
\endproof

\proclaim
Lemma 5.2.
Let $A$ be a commutative Hopf algebra over a perfect field $k$. There is a 
uniquely determined right coideal subalgebra $B$ of $A$ such that $AB^+$ 
coincides with the nil radical $\frn(A)$ of $A$. Moreover, $B$ is locally 
finite dimensional and $A$ is faithfully flat over $B$.
\endproclaim

\Proof.
Since the class of reduced commutative algebras over $k$ is closed under 
tensor products, the algebra $A/\frn(A)\ot A/\frn(A)$ is reduced. It follows 
that the nil radical of $A\ot A$ is equal to $\,\frn(A)\ot A+A\ot\frn(A)$, and 
therefore
$$
\De\bigl(\frn(A)\bigr)\sbs\frn(A)\ot A+A\ot\frn(A).
$$
Since the antipode $S$ of $A$ is an algebra antiautomorphism, it maps the nil 
radical $\frn(A)$ onto itself. Thus $\frn(A)$ is a Hopf ideal of $A$. The 
algebra $A$ represents an affine group scheme $G$, and the factor algebra 
$A/\frn(A)$ represents a closed group subscheme which coincides with the 
largest reduced subscheme $G\red$ of $G$.

For a right coideal subalgebra $B$ of $A$ we have $B^+\sbs\frn(A)$ if and only 
if all elements of $B$ are invariant with respect to the left coaction of 
$A/\frn(A)$ on $A$, i.e.,
$$
B\sbs\lco{A/\frn(A)}A=\{a\in A\mid\De(a)-1\ot a\in\frn(A)\ot A\}.
$$
Since the maximal ideal $B^+$ of such an algebra consists of nilpotent 
elements, every finitely generated subalgebra of $B$ is finite dimensional. A 
commutative Hopf algebra is flat over all its coideal subalgebras 
\cite{Ma-W94, Th. 3.4}. Since $B^+$ is the only prime ideal of $B$, it follows 
that $A$ is faithfully flat over $B$. Then $B=\lco{A/AB^+}A$ by Takeuchi's 
correspondence \cite{Tak79, Th.~1}. If $AB^+=\frn(A)$, then we must have 
$B=\lco{A/\frn(A)}A$.

It remains to prove that the equality $AB^+=\frn(A)$ does hold for the right 
coideal subalgebra $B=\lco{A/\frn(A)}A$. Since $A$ is a directed union of 
finitely generated Hopf subalgebras and $\frn(A')=\frn(A)\cap A'$ for each 
subalgebra $A'$, the proof reduces to the case when $A$ is finitely generated. 
In this case $G$ is a group scheme of finite type over $k$. By \cite{Dem-G, 
Ch. III, \S~3, Th. 5.4} there exists the quotient scheme $G\red\backslash G$, 
and $B$ is precisely the algebra of global sections of its structure sheaf. 
Since $G$ and $G\red$ have the same $K$-valued points for each extension field 
$K$ of the base field, the geometric realization of $G\red\backslash G$ 
consists of a single point. This implies that the scheme $G\red\backslash G$ is 
finite, and therefore affine. The equality $AB^+=\frn(A)$ expresses the fact 
that $G\red$ is the fiber of the canonical morphism $G\to G\red\backslash G$.
\endproof

\proclaim
Lemma 5.3.
Suppose that $B\sbs A\sbs H$ where $H$ is a weakly finite Hopf algebra, $A$ a 
commutative Hopf subalgebra, and $B$ a finite dimensional right coideal 
subalgebra whose ideal $B^+$ is nilpotent. Then a Hopf module $M\in\MAH$ will 
be projective in $\MA$ provided that $M/MB^+$ is projective in $\calM_{A/AB^+}$.
\endproclaim

\Proof.
By \cite{Sk07, Th. 6.1} all objects of the category $\MBH$ are free $B$-modules. 
In particular, this holds for all objects of $\MAH$. Consider any exact sequence
$$
0\to K\mapr{f}F\mapr{}M\to0
$$
in the category $\MA$ with a free $A$-module $F$. Since $M$ is $B$-free, this 
sequence splits in $\MB$. Making reduction modulo the ideal $B^+$, we obtain an 
exact sequence
$$
0\to K/KB^+\mapr{f'}F/FB^+\mapr{}M/MB^+\to0
$$
in $\calM_{A/AB^+}$. The assumption that $M/MB^+$ is projective in 
$\calM_{A/AB^+}$ implies that the latter sequence splits. Hence there exists 
an $A$-linear map $g:F\to K$ such that the induced map $F/FB^+\to K/KB^+$ 
gives a splitting of $f'$, i.e.,
$$
\Img(\,gf-\id_K)\sbs KB^+.
$$
Since $B^+$ is a nilpotent ideal, it follows that $gf-\id_K$ is a 
nilpotent endomorphism of the $A$-module $K$. But then $gf$ is an invertible 
endomorphism, and therefore
$$
F=\Img f\oplus\Ker g\,.
$$
Hence $M\cong\Ker g$ is a direct summand of a free $A$-module.
\endproof

\proclaim
Theorem 5.4.
Let $A$ be a commutative Hopf subalgebra of a weakly finite Hopf algebra $H$. 
Then all Hopf modules of the categories $\MAH$ and $\AMH$ are projective 
$A$-modules. As a consequence, $H$ is right and left faithfully flat over $A$. 
\endproclaim

\Proof.
Let $M\in\MAH$. Since $M$ is projective in $\MA$ if and only if $M\ot\kov$ is 
projective in $\calM_{A\ot\kov}$ where $\kov$ is the algebraic closure of the 
base field $k$, we may assume $k$ to be algebraically closed. As has been 
recalled in Lemma 1.1 it suffices to prove projectivity for $A$-finite Hopf 
modules.

Consider first the case when $A$ is finitely generated, and so the group 
scheme $G$ represented by $A$ has finite type over $k$. Let $B\sbs A$ be the 
right coideal subalgebra given by Lemma 5.2. Then $A/AB^+=A/\frn(A)$ is the 
algebra representing the largest reduced group subscheme $G\red$ which is 
smooth since $k$ is perfect. Hence $A/AB^+$ is a ring of finite global 
dimension. Since $B$ corresponds to the finite scheme $G\red\backslash G$, we 
have $\,\dim B<\infty$.

Let $M$ be $A$-finite. Since $A$ is Noetherian, $M$ is finitely presented in 
$\MA$. As in the proof of Theorem 4.4 take an infinite exact sequence
$$
\cdots\to F_2\to F_1\to F_0\to M\to0
$$
in $\MAH$ such that each $F_i$ is a finitely generated free $A$-module. Since 
all objects of the category $\MBH$ are free $B$-modules by \cite{Sk07, Th. 6.1}, 
this sequence splits in $\MB$. Hence the induced sequence of $A/AB^+$-modules
$$
\cdots\to F_2/F_2B^+\to F_1/F_1B^+\to F_0/F_0B^+\to M/MB^+\to0
$$
is exact as well. Now denote by $d$ the projective dimension of the 
$A/AB^+$-module $M/MB^+$. Suppose that $d>0$. Then the image $N$ of $F_d$ in 
$F_{d-1}$ is an $A$-finite $\MAH$-subobject of $F_{d-1}$. From the exact 
sequence $F_{d+1}\to F_d\to N\to0$ in $\MA$ we infer that $N/NB^+$ is 
isomorphic to the cokernel of the map 
$$
F_{d+1}/F_{d+1}B^+\to F_d/F_dB^+,
$$
which is projective in $\calM_{A/AB^+}$ by the definition of $d$. Lemma 5.3 
shows that $N$ is projective in $\MA$. Hence $N$ is an $\MA$-direct summand of 
$F_{d-1}$ by Lemma 5.1, and it follows that $F_{d-1}/N$ is projective in $\MA$ 
as well. This means that $\,\pd_AM<d$, but the projective dimension of the 
$A/AB^+$-module $M/MB^+$ cannot exceed that of the $A$-module $M$, and we 
arrive at a contradiction. Thus $d=0$, and we deduce that $M$ is projective in 
$\MA$ applying Lemma 5.3 again.

In the case when $A$ is not finitely generated we conclude that each Hopf 
module $M\in\MAH$ is a flat $A$-module noting that the functor $M\ot_A?$ is 
naturally isomorphic to the directed colimit $\,\limdir\,(M\ot_{A'}?)$ over 
the set of finitely generated Hopf subalgebras $A'\sbs A$. Hence $M$ is 
projective in $\MA$ whenever $M$ is finitely presented in $\MA$. Arguing 
as in the proof of Theorem 4.4 we deduce that each $A$-finite object of 
$\MAH$ is finitely presented in $\MA$, and this implies that all objects of 
$\MAH$ are projective in $\MA$. Finally, the projectivity conclusion is 
extended to the category $\AMH$ as in the proof of Corollary 4.7.
\endproof

\section
Appendix. Topological lemmas on countable irreducibility

A topological space $X$ is a \emph{Baire space} if every countable 
intersection of dense open subsets is dense in $X$. In an \emph{irreducible 
space} all nonempty open subsets are dense. Therefore a topological space $X$ 
is an irreducible Baire space if and only if every countable intersection of 
nonempty open subsets is nonempty. Equivalently, $X$ cannot be presented as 
the union of any countable collection of proper closed subsets. We call such 
spaces \emph{countably irreducible}.

Examples of such spaces arise naturally in algebraic geometry. For this we 
have to assume everywhere that
$$
\text{\emph{the base field $\,k\,$ is uncountable and algebraically closed.}}
$$
We will view each algebraic variety over the base field $k$ as a set of closed 
points equipped with Zariski topology.

\proclaim
Lemma A1.
Any irreducible algebraic variety $X$ with its Zariski topology is countably 
irreducible.
\endproclaim

\Proof.
By passing to an affine open subset of $X$ the proof is reduced to the case 
where $X$ is affine. By Noether's normalization lemma there exists then a 
finite surjective morphism $f:X\to\bbA^n$ where $\bbA^n$ is the affine space 
of dimension $n=\dim X$. Suppose that $X$ is the union of a sequence ot its 
closed subsets $Z_1,Z_2,\ldots$. Then $\bbA^n=\bigcup_{i=1}^\infty f(Z_i)$, 
and each subset $f(Z_i)$ is closed in $\bbA^n$. Since $\,\dim f(Z_i)=\dim Z_i$, 
we have $f(Z_i)\ne\bbA^n$ whenever $Z_i\ne X$.

Thus it suffices to prove that $\bbA^n$ is countably irreducible. We proceed 
by induction on $n$. Suppose that $\bbA^n=\bigcup_{i=1}^\infty W_i$ where 
$W_1,W_2,\ldots$ are proper closed subsets of $\bbA^n$. By Hilbert's 
Nullstellensatz we can find for each $i$ a nonzero polynomial function 
$f_i:\bbA^n\to k$ vanishing on $W_i$. The zero set $V(f_i)$ of $f_i$ is a 
closed subvariety in $\bbA^n$. It has finitely many irreducible components. 
Now let $\ell:\bbA^n\to k$ be any nonzero linear function, and for each 
$\al\in k$ put 
$$
Y_\al=\{x\in\bbA^n\mid\ell(x)=\al\}.
$$
Clearly $Y_\al$ is an affine subspace of dimension $n-1$, and $Y_\al\ne Y_\be$ 
whenever $\al\ne\be$. Since $k$ is uncountable, there exists $\al$ such that 
$Y_\al$ coincides with none of the irreducible components of the countable 
collection of hypersurfaces $V(f_i)$. Then $Y_\al\not\sbs W_i$ for each $i$, 
whence $Y_\al=\bigcup_{i=1}^\infty Y_\al\cap W_i$ is a countable union of 
proper closed subsets of $Y_\al$. Since $Y_\al\cong\bbA^{n-1}$, the induction 
argument completes the proof.
\endproof

Suppose now that $\,X=\prod_{i\in I}X_i\,$ is the cartesian product of an 
indexed collection $(X_i)_{i\in I}$ of algebraic varieties. For 
each subset $E\sbs I$ put
$$
X_E=\prod_{i\in E}X_i
$$
and denote by $p_E$ the projection $X\to X_E$. If $F\sbs I$ is a finite 
subset, then $X_F$ has the structure of an algebraic variety which makes $X_F$ 
the direct product of varieties $X_i$, $i\in F$, in the category of algebraic 
varieties.

We equip $X$ with the coarsest topology such that $\,p_F:X\to X_F\,$ is 
continuous for each finite subset $F\sbs I$. Its base of open sets consists of 
all subsets of $X$ of the form $p_F^{-1}(O)$ where $F$ is a finite subset of 
$I$ and $O$ is a Zariski open subset of $X_F$. In other words, this topology 
makes $X$ homeomorphic to the inverse limit $\,\liminv_FX_F\,$ in the category 
of topological spaces. This topology on $X$ which we call the 
\emph{pro-Zariski topology} is finer than the direct product topology.

Note that for any subset $E\sbs I$ the projection $p_E:X\to X_E$ is an open 
map, i.e., $p_E$ maps open subsets of $X$ onto open subsets of $X_E$. If 
$U=p_F^{-1}(O)$ where $F$ is a finite subset of $I$ and $O$ is an open subset 
of $X_F$, then the openness of $p_E(U)$ in $X_E$ can be seen by considering 
the commutative diagram
$$
\diagram{
X & \lmapr2{p_F} & X_F \cr
\noalign{\smallskip}
\mapd{p_E}{} && \mapd{}{g} \cr
\noalign{\smallskip}
X_E & \lmapr2{h} & X_{E\cap F} \cr
}
$$
in which all maps are the canonical projections. The projection 
$g:X_F\to X_{E\cap F}$ is a flat morphism of algebraic varieties, and 
therefore it is open. Hence $g(O)$ is open in $X_{E\cap F}$, and 
$p_E(U)=h^{-1}\bigl(g(O)\bigr)$ is indeed open in $X_E$.

\proclaim
Lemma A2.
The cartesian product $X=\prod_{i\in I}X_i$ of any family of irreducible 
algebraic varieties is countably irreducible in the pro-Zariski topology.
\endproclaim

\Proof.
We have to prove that $\bigcap_{n=1}^\infty D_n\ne\varnothing$ for any 
sequence $D_1,D_2,\ldots$ of nonempty open subsets of $X$. It suffices to do 
this assuming that each $D_n$ is in the base of topology, i.e., 
$D_n=p_{F_n}^{-1}(O_n)$ where $F_n$ is a finite subset of $I$ and $O_n$ is 
a nonempty open subset of $X_{F_n}$. 

If $E\sbs I$ is any subset, then $p_E(D_n)$ is a nonempty open subset of $X_E$ 
since $p_E$ is an open map. Denote by $K$ the set of all pairs $(E,x)$ where
$$
E\sbs I\quad\text{and}\quad x=(x_i)_{i\in E}\in X_E
$$
are such that $x\in\bigcap_{n=1}^\infty p_E(D_n)$. Define a partial order on 
$K$ setting $(E,x)\le(E',x')$ if $E\sbs E'$ and $x_i=x'_i$ for all $i\in E$, 
i.e., $x=x'_E$ where $x'_E$ denotes the projection of $x'$ to $X_E$. We 
claim that $K$ satisfies the hypotheses of Zorn's lemma. If $E$ is finite, 
then $X_E$ is an irreducible algebraic variety, whence $\bigcap_{n=1}^\infty 
p_E(D_n)\ne\varnothing$ by Lemma A1. This shows that $K\ne\varnothing$.

Suppose that $C$ is some chain in $K$. Denote by $J$ the union of all subsets 
$E\sbs I$ which occur as the first components of pairs $(E,x)\in C$. If $i\in J$, 
then there is $y_i\in X_i$ such that $y_i=x_i$ for each pair $(E,x)\in C$ with 
$i\in E$. This gives an element $y=(y_i)_{i\in J}\in X_J$ such that for each 
$(E,x)\in C$ we have $E\sbs J$ and $x=y_E$. Let us show that $y\in p_J(D_n)$. 
Note that this containment means precisely that there exists 
$z\in O_n\sbs F_n$ such that $y_i=z_i$ for all $i\in J\cap F_n$. Since $J\cap 
F_n$ is a finite subset of $J$, there exists some pair $(E,x)\in C$ such that 
$J\cap F_n\sbs E$. But $x\in p_E(D_n)$ by the definition of $K$, and therefore 
there exists $z\in O_n$ such that $x_i=z_i$ for all 
$i\in E\cap F_n=J\cap F_n$. Since $y_i=x_i$ for all $i\in E$, we do get the 
required element $z$. Hence $y\in\bigcap_{n=1}^\infty p_J(D_n)$, and therefore 
$(J,y)\in K$. It is now clear that $(J,y)$ is the supremum of $C$ in $K$.

Thus Zorn's lemma does apply. It shows that $K$ has a maximal element. Let now 
$(E,x)$ be any maximal element of $K$, and suppose that $E\ne I$. We may 
identify $X$ with $X_E\times X_{I\setm E}$ set-theoretically. Since the map 
$X_{I\setm E}\to X$, $\,t\mapsto(x,t)$, is continuous, the subset 
$$
D'_n=\{t\in X_{I\setm E}\mid(x,t)\in D_n\}
$$
is open in $X_{I\setm E}$ for each $n$. Moreover, $D'_n\ne\varnothing$ since 
$x\in p_E(D_n)$. Denoting by $p'_i$ the projection $X_{I\setm E}\to X_i$ for 
$i\in I\setm E$ we have $\bigcap_{n=1}^\infty p'_i(D'_n)\ne\varnothing$ by 
Lemma A1. Taking any element $u\in X_i$ lying in this intersection and setting 
$$
E'=E\cup\{i\},\quad x'=(x,u)\in X_{E'},
$$
we get $x'\in\bigcap_{n=1}^\infty p_{E'}(D_n)$. Hence $(x',E')\in K$, but this 
contradicts maximality of the pair $(x,E)$ in $K$.

It follows that $K$ contains a pair $(E,x)$ with $E=I$, and then 
$x\in\bigcap_{n=1}^\infty D_n$.
\endproof

\proclaim
Lemma A3.
Let $X=\prod_{i\in I}X_i$ and $Y=\prod_{j\in J}Y_j$ where $(X_i)_{i\in I}$ and 
$(Y_j)_{j\in J}$ are two collections of algebraic varieties. 
Suppose that $\ph:X\to Y$ is a map such that for each $j\in J$ there exist a 
finite subset $F_j\sbs I$ and a morphism of algebraic varieties 
$\,\ph_j:X_{F_j}\to Y_j\,$ rendering commutative the diagram
$$
\diagram{
X & \lmapr2{p_{F_j}} & X_{F_j} \cr
\noalign{\smallskip}
\mapd{\ph}{} && \mapd{}{\ph_j} \cr
\noalign{\smallskip}
Y & \lmapr2{q_j} & Y_j \cr
}
$$
where $p_{F_j}$ and $q_j$ are the canonical projections. Then $\ph$ is 
continuous.
\endproclaim

\Proof.
It suffices to prove that $\ph^{-1}(U)$ is open in $X$ for each open set $U$ 
in the standard base of topology on $Y$.

So let $U=q_{F'}^{-1}(O)$ where 
$F'\sbs J$ is a finite subset, $q_{F'}:Y\to Y_{F'}$ the projection, and 
$O\sbs Y_{F'}$ an open subset. Consider the finite subset 
$F=\bigcup_{j\in F'}F_j$ of $I$ and the map $\ph_F:X_F\to Y_{F'}$ such that 
for each $j\in F'$ there is a commutative diagram
$$
\diagram{
X & \lmapr2{p_F} & X_F & \lmapr2{} & X_{F_j} \cr
\noalign{\smallskip}
\mapd{\ph}{} && \mapd{}{\ph_F} && \mapd{}{\ph_j} \cr
\noalign{\smallskip}
Y & \lmapr2{q_{F'}} & Y_{F'} & \lmapr2{} & Y_j \cr
}
$$
where all horizontal arrows are projections. In fact there is a unique map 
$\ph_F$ making commutative the right square in the diagram for each $j\in F'$, 
and then the left square commutes as well. Since $Y_{F'}=\prod_{j\in F'}Y_j$ 
is a direct product in the category of algebraic varieties, the map $\ph_F$ is 
a morphism of algebraic varieties. In particular, $\ph_F$ is continuous. Hence 
$\ph_F^{-1}(O)$ is open in $X_F$, and 
$\ph^{-1}(U)=p_F^{-1}\bigl(\ph_F^{-1}(O)\bigr)$ is open in $X$.
\endproof

\proclaim
Lemma A4.
Let $X=\prod_{i\in I}X_i$ be the cartesian product of a collection 
$(X_i)_{i\in I}$ of finite-dimensional vector spaces, and let $L$ be a vector 
subspace of finite codimension in $X$. Then $L$ is closed in $X$ if and only 
if $\,\Ker p_F\sbs L\,$ for some finite subset $F\sbs I$.
\endproclaim

\Proof.
Every closed subset of $X$ is the intersection of a collection of subsets 
of the form $p_F^{-1}(Z)$ where $F$ is a finite subset of $I$ and $Z$ is a 
Zariski closed subset of the finite-dimensional vector space $X_F$. Note that 
$L\sbs p_F^{-1}(Z)$ if and only if $p_F(L)\sbs Z$, if and only if 
$\,p_F^{-1}\bigl(p_F(L)\bigr)\sbs p_F^{-1}(Z)$.

Therefore $L$ is closed if and only if $\,L=\bigcap\,p_F^{-1}\bigl(p_F(L)\bigr)$, 
the intersection over all finite subsets of $I$. For each $F$ the inverse 
image $p_F^{-1}\bigl(p_F(L)\bigr)$ is a vector subspace of $X$ containing $L$. 
Since $L$ has finite codimension in $X$ and the map 
$F\mapsto p_F^{-1}\bigl(p_F(L)\bigr)$ reverses inclusions, there exists a 
smallest subspace among all subspaces $p_F^{-1}\bigl(p_F(L)\bigr)$ associated 
with finite subsets  of $I$.

Hence $L$ is closed in $X$ if and only if $\,L=p_F^{-1}\bigl(p_F(L)\bigr)\,$ 
for some finite subset $F$ of the index set $I$. It remains to note that the 
previous equality is equivalent to the inclusion $\,\Ker p_F\sbs L$.
\endproof

\setitemsize(a)
\proclaim
Lemma A5.
Let $X=\prod_{i\in I}X_i$ and $Y=\prod_{j\in J}Y_j$ where $(X_i)_{i\in I}$ and 
$(Y_j)_{j\in J}$ are two collections of finite-dimensional vector spaces. 
For a linear map $\ph:X\to Y$ the following conditions are equivalent:

\item(a)
$\ph$ is continuous,

\item(b)
$\Ker(q_j\circ\ph)$ is a closed subspace of $X$ for each $j\in J,$

\item(c)
for each $j\in J$ there exists a finite subset $F\sbs I$ such that 
$$
\Ker p_F\sbs\Ker(q_j\circ\ph)
$$
where $\,p_F:X\to X_F\,$ and $\,q_j:Y\to Y_j\,$ are the canonical projections.
\endproclaim

\Proof.
Since all single element subsets of $Y_j$ are closed in the Zariski topology 
and $q_j$ is continuous, the kernel of $q_j$ is closed in  $Y$. If $\ph$ is 
continuous, then
$$
\Ker(q_j\circ\ph)=\ph^{-1}(\Ker q_j)
$$
is closed in  $X$. Hence (a)$\,\Rar\,$(b). By Lemma A4 (b)$\,\Rar\,$(c) since 
$\Ker(q_j\circ\ph)$ is a vector subspace of finite codimension in $X$.

Suppose that (c) holds. Then the linear map $q_j\circ\ph:X\to Y_j$ factors as 
$\psi\circ p_F$ where $\psi:X_F\to Y_j$ is a linear map. Hence (a) follows 
from Lemma A3.
\endproof

\proclaim
Corollary A6.
Let $X$ and $Y$ be as in Lemma A5, and let $\ph:X\to Y$ be an invertible 
linear map. Then $\ph$ is bicontinuous, i.e., both $\ph$ and its inverse 
$\ph^{-1}$ are continuous, if and only if the assignment $\,L\mapsto\ph(L)\,$ 
gives a bijection between the sets of closed vector subspaces of finite 
codimension in $X$ and in $Y$.
\endproclaim

\references
\nextref Ag11
\auth{A.L.,Agore}
\paper{Categorical constructions for Hopf algebras}
\journal{Comm. Algebra}
\Vol{39}
\Year{2011}
\Pages{1476-1481}

\nextref Ar-G03
\auth{S.,Arkhipov;D.,Gaitsgory}
\paper{Another realization of the category of modules over the small quantum group}
\journal{Adv. Math.}
\Vol{173}
\Year{2003}
\Pages{114-143}

\nextref Ber74
\auth{G.M.,Bergman}
\paper{Modules over coproducts of rings}
\journal{Trans. Amer. Math. Soc.}
\Vol{200}
\Year{1974}
\Pages{1-32}

\nextref Br98
\auth{K.A.,Brown}
\paper{Representation theory of Noetherian Hopf algebras satisfying a polynomial identity}
\InBook{Trends in the representation theory of finite dimensional algebras}
\BkSer{Contemp. Math.}
\BkVol{229}
\publisher{Amer. Math. Soc.}
\Year{1998}
\Pages{49-79}

\nextref Br-G97
\auth{K.A.,Brown;K.R.,Goodearl}
\paper{Homological aspects of Noetherian PI Hopf algebras and irreducible modules of maximal dimension}
\journal{J.~Algebra}
\Vol{198}
\Year{1997}
\Pages{240-265}

\nextref Br-Z10
\auth{K.A.,Brown;J.J.,Zhang}
\paper{Prime regular Hopf algebras of GK-dimension one}
\journal{Proc. London Math. Soc.}
\Vol{101}
\Year{2010}
\Pages{260-302}

\nextref Cae-G04
\auth{S.,Caenepeel;T.,Gu\'ed\'enon}
\paper{Projectivity of a relative Hopf module over the subring of coinvariants}
\InBook{Hopf algebras}
\publisher{Marcel Dekker}
\Year{2004}
\Pages{97-108}

\nextref Dem-G
\auth{M.,Demazure;P.,Gabriel}
\book{Groupes Alg\'ebriques I}
\publisher{Masson}
\Year{1970}

\nextref Doi83
\auth{Y.,Doi}
\paper{On the structure of relative Hopf modules}
\journal{Comm. Algebra}
\Vol{11}
\Year{1983}
\Pages{243-255}

\nextref Doi92
\auth{Y.,Doi}
\paper{Unifying Hopf modules}
\journal{J.~Algebra}
\Vol{153}
\Year{1992}
\Pages{373-385}

\nextref Doi-T86
\auth{Y.,Doi;M.,Takeuchi}
\paper{Cleft comodule algebras for a bialgebra}
\journal{Comm. Algebra}
\Vol{14}
\Year{1986}
\Pages{801-817}

\nextref Goo-W
\auth{K.R.,Goodearl;R.B.,Warfield Jr.}
\book{An Introduction to Noncommutative Noetherian Rings}
second edition,
\publisher{Cambridge Univ. Press}
\Year{2004}

\nextref Goo-Z10
\auth{K.R.,Goodearl;J.J.,Zhang}
\paper{Noetherian Hopf algebra domains of Gelfand-Kirillov dimension two}
\journal{J.~Algebra}
\Vol{324}
\Year{2010}
\Pages{3131-3168}

\nextref Lar-Sw69
\auth{R.G.,Larson;M.E.,Sweedler}
\paper{An associative orthogonal bilinear form for Hopf algebras}
\journal{Amer. J. Math.}
\Vol{91}
\Year{1969}
\Pages{75-94}

\nextref Ma94
\auth{A.,Masuoka}
\paper{Quotient theory of Hopf algebras}
\InBook{Advances in Hopf algebras}
\publisher{Marcel Dekker}
\Year{1994}
\Pages{107-133}

\nextref Ma-W94
\auth{A.,Masuoka;D.,Wigner}
\paper{Faithful flatness of Hopf algebras}
\journal{J.~Algebra}
\Vol{170}
\Year{1994}
\Pages{156-164}

\nextref Mo
\auth{S.,Montgomery}
\book{Hopf Algebras and their Actions on Rings}
\publisher{Amer. Math. Soc.}
\Year{1993}

\nextref Nich-Z89a
\auth{W.D.,Nichols;M.B.,Zoeller}
\paper{A Hopf algebra freeness theorem}
\journal{Amer. J. Math.}
\Vol{111}
\Year{1989}
\Pages{381-385}

\nextref Nich-Z89b
\auth{W.D.,Nichols;M.B.,Zoeller}
\paper{Freeness of infinite dimensional Hopf algebras over grouplike subalgebras}
\journal{Comm. Algebra}
\Vol{17}
\Year{1989}
\Pages{413-424}

\nextref Ox61
\auth{J.C.,Oxtoby}
\paper{Cartesian products of Baire spaces}
\journal{Fund. Math.}
\Vol{49}
\Year{1961}
\Pages{157-166}

\nextref Par72
\auth{B.,Pareigis}
\paper{On the cohomology of modules over Hopf algebras}
\journal{J.~Algebra}
\Vol{22}
\Year{1972}
\Pages{161-182}

\nextref Por11
\auth{H.-E.,Porst}
\paper{Limits and colimits of Hopf algebras}
\journal{J.~Algebra}
\Vol{328}
\Year{2011}
\Pages{254-267}

\nextref Rot
\auth{J.J.,Rotman}
\book{An Introduction to Homological Algebra}
second edition,
\publisher{Springer}
\Year{2009}

\nextref Scha00
\auth{P.,Schauenburg}
\paper{Faithful flatness over Hopf subalgebras: counterexamples}
\InBook{Interactions between ring theory and representations of algebras}
\publisher{Marcel Dekker}
\Year{2000}
\Pages{331-344}

\nextref Schn92
\auth{H.-J.,Schneider}
\paper{Normal basis and transitivity of crossed products for Hopf algebras}
\journal{J.~Algebra}
\Vol{152}
\Year{1992}
\Pages{289-312}

\nextref Schn93
\auth{H.-J.,Schneider}
\paper{Some remarks on exact sequences of quantum groups}
\journal{Comm. Algebra}
\Vol{21}
\Year{1993}
\Pages{3337-3357}

\nextref Sk06
\auth{S.,Skryabin}
\paper{New results on the bijectivity of antipode of a Hopf algebra}
\journal{J.~Algebra}
\Vol{306}
\Year{2006}
\Pages{622-633}

\nextref Sk07
\auth{S.,Skryabin}
\paper{Projectivity and freeness over comodule algebras}
\journal{Trans. Amer. Math. Soc.}
\Vol{359}
\Year{2007}
\Pages{2597-2623}

\nextref Sk08
\auth{S.,Skryabin}
\paper{Projectivity of Hopf algebras over subalgebras with semilocal central localizations}
\journal{J.~$K$-Theory}
\Vol{2}
\Year{2008}
\Pages{1-40}

\nextref Sk21
\auth{S.,Skryabin}
\paper{Flatness of Noetherian Hopf algebras over coideal subalgebras}
\journal{Algebr. Represent. Theory}
\Vol{24}
\Year{2021}
\Pages{851-875}

\nextref Sk-
\auth{S.,Skryabin}
\paper{On Takeuchi's correspondence}
arXiv: 2501.06045

\nextref Tak71
\auth{M.,Takeuchi}
\paper{Free Hopf algebras generated by coalgebras}
\journal{J.~Math. Soc. Japan}
\Vol{23}
\Year{1971}
\Pages{561-582}

\nextref Tak72
\auth{M.,Takeuchi}
\paper{A correspondence between Hopf ideals and sub-Hopf algebras}
\journal{Manuscripta Math.}
\Vol{7}
\Year{1972}
\Pages{251-270}

\nextref Tak79
\auth{M.,Takeuchi}
\paper{Relative Hopf modules---equivalences and freeness criteria}
\journal{J.~Algebra}
\Vol{60}
\Year{1979}
\Pages{452-471}

\nextref Wu-Z02
\auth{Q.-S.,Wu;J.J.,Zhang}
\paper{Regularity of involutory PI Hopf algebras}
\journal{J.~Algebra}
\Vol{256}
\Year{2002}
\Pages{599-610}

\endreferences
\bye